\documentclass[a4j,10pt]{article}
\usepackage{amsmath}
\usepackage{amssymb}
\usepackage{amsthm}
\usepackage{graphics}
\usepackage[dvips]{graphicx}
\usepackage{amscd}

\numberwithin{equation}{section}
\theoremstyle{definition}
\newtheorem{example}{example}[section]
\newtheorem{df}[example]{Definition}
\theoremstyle{remark}
\newtheorem{rem}[example]{
{\bf Remark}}
\theoremstyle{plain}

\newtheorem{prop}[example]{Proposition}

\newtheorem{theorem}{Theorem}[section]

\newtheorem{lemma}[theorem]{Lemma}

\newtheorem{remark}[theorem]{Remark}

\theoremstyle{remark}

\begin{document}

\title{Algebraic Ricci Solitons of three-dimensional Lorentzian Lie groups}
\author{ W. Batat and K. Onda\thanks{%
First author was supported by funds of JSPS and ENSET d'Oran. The second
author was supported by funds of Nagoya University and OCAMI. \newline
} }
\maketitle

\begin{abstract}
We classify Algebraic Ricci Solitons of three-dimensional Lorentzian Lie
groups. All algebraic Ricci solitons that we obtain are solvsolitons. In
particular, we obtain new solitons on $G_{2}$, $G_{5}$, and $G_{6}$, and we
prove that, contrary to the Riemannian case, Lorentzian Ricci solitons need
not to be algebraic Ricci solitons. 
\end{abstract}










\section{\protect\LARGE Introduction and preliminaries }

\setcounter{equation}{0}

The concept of the algebraic Ricci soliton was first introduced by Lauret in
Riemannian case (\cite{L01}). The definition extends to the
pseudo-Riemannian case:

\begin{df}
\textit{Let } $\left( \mathit{G,g}\right) $\textit{\ be a simply connected
Lie group equipped with the left-invariant pseudo-Riemannian metric $g$, and
let } ${\mathfrak{g}}$ denote \textit{the Lie algebra }of $\mathit{G.}$ 
\textit{Then $g$ is called an algebraic Ricci soliton if it satisfies } 
\begin{equation}
\mathrm{Ric}=c\mathrm{Id}+D  \label{Soliton}
\end{equation}%
\textit{where }$\mathrm{Ric}$\textit{\ denotes the Ricci operator, $c$ is a
real number, and }$\mathrm{D\in Der}\left( \mathfrak{g}\right) $ ($\mathrm{D}
$ is \textit{a derivation of }$\mathit{\mathfrak{g}}$\textit{)}, that is:%
\begin{equation*}
\mathit{\mathrm{D}[X,Y]=[\mathrm{D}X,Y]+[X,\mathrm{D}Y]\ \,}\text{\ \textit{%
for any }}\mathit{\,\ X,Y\in \mathfrak{g}.}
\end{equation*}%
\textit{\ In particular, an algebraic Ricci soliton on a solvable Lie group,
(a nilpotent Lie group) is called a solvsoliton (a nilsoliton).}
\end{df}

Obviously, Einstein metrics on a Lie group are algebraic Ricci solitons.
From now on, by a trivial algebraic Ricci soliton we shall mean an algebraic Ricci soliton which
is Einstein. 

\begin{rem}
Let $G$ be a semi-simple Lie group, $g$ a left-invariant Riemannian metric.
If $g$ is an algebraic Ricci soliton, then $g$ is Einstein (\cite{L01}).
\end{rem}

Next we introduce Ricci solitons. Let $g_{0}$ be a pseudo-Riemannian metric
on a manifold $M^{n}$. If $g_{0}$ satisfies 
\begin{equation}
\varrho \lbrack g_{0}]=cg_{0}+L_{X}g_{0}\ ,  \label{Ricci soliton}
\end{equation}%
where $\varrho $ denotes the Ricci tensor, $X$ is a vector field and $c$ is
a real constant, then $(M^{n},g_{0},X,c)$ is called a \textit{Ricci soliton
structure} and $g_{0}$ is \textit{the Ricci soliton.} Moreover we say that
the Ricci soliton $g_{0}$ is a \textit{gradient Ricci soliton} if the vector
field $X$ satisfies $X=\nabla f$, where $f$ is a function. The Ricci soliton 
$g_{0}$ is said to be a \textit{non-gradient Ricci soliton} if the vector
field $X$ satisfies $X\neq \nabla f$ for any function $f$. If $c$ is
positive, zero, or negative, then $g_{0}$ is called a shrinking, steady, or
expanding Ricci soliton, respectively. According to \cite{CK04}, we check
that a Ricci soliton is a Ricci flow solution.

\begin{prop}[\protect\cite{CK04}]
\textit{A pseudo-Riemannian metric $g_{0}$ is Ricci soliton if and only if $%
g_{0}$ is the initial metric of the Ricci flow equation,} 
\begin{equation*}
\frac{\partial }{\partial t}g(t)_{ij}=-2 \varrho \lbrack g(t)]_{ij}\ ,
\end{equation*}%
\textit{and the solution is expressed as $g(t)=c(t)(\varphi _{t})^{\ast
}g_{0}$, where $c(t)$ is a scaling parameter, and $\varphi _{t}$ is a
diffeomorphism.}
\end{prop}

An interesting example of Ricci solitons is $(\mathbb{R}^{2},g_{st}\ ,X,c),$
where the metric $g_{st}$ is the Euclidean metric on $\mathbb{R}^{2}$, the
vector field $X$ is 
\begin{equation*}
X=\nabla f,\quad f=\dfrac{\left\vert x\right\vert ^{2}}{2}
\end{equation*}%
and $c$ is a real number. This is a gradient Ricci soliton structure and so, 
$g_{st}$ is the gradient Ricci soliton, named Gaussian soliton. 
In the closed Riemannian case, Perelman \cite{P02} proved that any Ricci
soliton is a gradient Ricci soliton, and any steady or expanding Ricci
soliton is an Einstein metric with the Einstein constant zero or negative,
respectively. However in the non-compact Riemannian case, a Ricci soliton is
not necessarily gradient and a steady or expanding Ricci soliton is not
necessarily Einstein. In fact, any left-invariant Riemannian metric on the
three-dimensional Heisenberg group is an expanding non-gradient Ricci
soliton which is not an Einstein metric (see \cite{BD07}, \cite{GIK06}, \cite%
{L07}). In the Riemannian case, all homogeneous non-trivial Ricci solitons
are expanding Ricci solitons. In the pseudo-Riemannian case, there are
shrinking homogeneous non-trivial Ricci solitons discovered in \cite{O10},
the vectors fields of these Ricci solitons are not left-invariant.

In \cite{L01}, Lauret studied the relation between algebraic Ricci solitons and Ricci
solitons on Riemannian manifolds. More precisely, he proved that any
left-invariant Riemannian algebraic Ricci soliton metric is a Ricci soliton. This was
extended by the second author to the pseudo-Riemannian case :

\begin{theorem}[\protect\cite{L01}, \protect\cite{O11}]\label{O11}
\label{RS} Let $\left( G,g\right) $ be a simply connected Lie group endowed
with a left-invariant pseudo-Riemannian metric $g$. If $g$ is an algebraic Ricci soliton,
then $g$ is the Ricci soliton, that is, $g$ satisfy 
\eqref{Ricci soliton}%
, such that%
\begin{equation*}
X=\dfrac{d\varphi _{t}}{dt}\Big|_{t=0}(p)\text{ \ and }\mathrm{\exp }\left( 
\frac{t}{2}\mathrm{D}\right) =d\varphi _{t}|_{e},
\end{equation*}%
where $e$ denotes the identity element of $G$ .
\end{theorem}


On the other hand, it was proven in \cite{Di} that three-dimensional Lie
groups do not admit left-invariant Riemannian Ricci solitons, where a
left-invariant Ricci soliton means that the metric and the vector field are
left-invariant. Recently, M. Brozos-Vazquez et. al studied the corresponding
existence problem in Lorentzian signature \cite{BCGG09}, and proving the
existence of expanding, steady and shrinking left-invariant Ricci solitons
on three-dimensional Lorentzian homogeneous manifolds. Furthermore, other
results of Lorentzian Ricci solitons are found in \cite{BBGG10}, \cite{BGG11}%
. 
In this paper, we shall give the complete classification of algebraic Ricci
solitons in three-dimensional unimodular Lie groups and in non-symmetric
non-unimodular Lie groups. 
We prove

\begin{theorem}
\label{maintheorem1} Let $(G,g)$ be a three-dimensional connected Lie group,
equipped with a left-invariant Lorentzian metric $g.$ The following are all
non-trivial Lorentzian algebraic Ricci solitons :

\begin{itemize}
\item $G=G_{1}$ with $\beta =0$.

\item $G=G_{2}$ with $\alpha =\beta =0$.

\item $G=G_{3}$ with

\begin{itemize}
\item $\beta =\gamma =0$ and $\alpha >0$, or $\gamma <0$ and $\alpha =\beta
=0$, or

\item $\gamma =0$ and $\beta =-\alpha <0$, or

\item $\beta =0$ and $\gamma =-\alpha <0$.
\end{itemize}

\item $G=G_{5}$ with

\begin{itemize}
\item $\left( \beta ,\gamma \right) \neq \left( 0,0\right) $ and $\alpha
^{2}+\beta ^{2}=\gamma ^{2}+\delta ^{2}$, or

\item $\left( \beta ,\gamma \right) =\left( 0,0\right) $.
\end{itemize}

\item $G=G_{6}$ with

\begin{itemize}
\item $\left( \beta ,\gamma \right) \neq \left( 0,0\right) $ and $\alpha
^{2}-\beta ^{2}=\delta ^{2}-\gamma ^{2}$, or

\item $\left( \beta ,\gamma \right) =\left( 0,0\right) $.
\end{itemize}

\item $G=G_{7}$ with $\gamma =0$.
\end{itemize}
\end{theorem}

It is natural to consider the opposite of Theorem \ref{O11}. 
A Riemannian manifold $(M, g)$ is called a solvmanifold if there exists a transitive solvable group of isometries.
Jablonski proved 
\begin{theorem}[\cite{J11}] 
Consider a solvmanifold $(M, g)$ which is a Ricci soliton. 
Then $(M, g)$ is isometric to a solvsoliton and the transitive solvable group may be chosen to be completely solvable. 
\end{theorem} 

In the Lorentzian case, there exists a Ricci soliton which is not an algebraic Ricci soliton.
Indeed, from the Remark $\ref{maintheorem2a}$ and $\ref{maintheorem2b}$, we obtain 

\begin{theorem}
There exist Lorentzian Ricci solitons which are not algebraic Ricci solitons on $SL(2,R)$.
\end{theorem}

The paper is organized in the following way. In Section 2 we collect some
basic facts concerning three-dimensional Lie groups equipped with a
left-invariant Lorentzian metrics. Algebraic Ricci solitons of Lorentzian Lie
algebras $\mathfrak{g}_{1}\mathfrak{-g}_{7}$ listed in Theorem \ref%
{unimodular} and Theorem \ref{non-unimodular}, will be then investigated in
Section 3 and 4. 
In particular, we obtain a new solitons on $G_{2}$, $G_{5}$ and $G_{6}$ and
we prove that there exist Ricci solitons but not algebraic Ricci solitons on 
$G_{1}=SL(2,\mathbb{R})$ and $G_{4}=SL(2,\mathbb{R})$. 

\section{${\protect\LARGE 3}${\protect\LARGE -dimensional Lorentzian Lie
groups}}

\setcounter{equation}{0}

A Lie group $G$ is said to be \textit{unimodular} if its left-invariant Haar
measure is also right-invariant. Milnor \cite{M76} gave an infinitesimal
reformulation of unimodularity for three-dimensional Lie groups, obtaining a
complete classification of three-dimensional unimodular Lie groups equipped
with a left-invariant Riemannian metric. Three-dimensional Lorentzian Lie
groups have been studied by \cite{R92}, S. Rahmani (see also \cite{CP}), who
introduced a cross product $X\times Y$ adapted to the Lorentzian case,
proving that the bracket operation on a three-dimensional Lie algebra is
related to the Lorentzian cross product by 
\begin{equation*}
\lbrack X,Y]=L(X\times \ Y),
\end{equation*}%
for a unique linear endomorphism $L$ on the Lie algebra. Further, the
three-dimensional Lie algebra is unimodular if and only if $L$ is
selfadjoint (\cite{R92} and \cite{CP}). This characterization together with
the canonical forms for Lorentzian selfadjoint operators given in \cite{O}
(pp. 261-262, ex. 19), leads to the classification of three-dimensional
Lorentzian unimodular Lie groups \cite{R92}, \cite{CP}:

\begin{theorem}
\label{unimodular} Let $(G,g)$ be a three-dimensional connected unimodular
Lie group, equipped with a left-invariant Lorentzian metric. Then, there
exists a pseudo-orthonormal frame field $\{e_{1},e_{2},e_{3}\}$, with $e_{3}$
timelike, such that the Lie algebra of $G$ is one of the following:

- $\left( \mathfrak{g}_{1}\right) :$%
\begin{eqnarray}
\lbrack e_{1},e_{2}] &=&\alpha e_{1}-\beta e_{3},  \label{g1} \\
\lbrack e_{1},e_{3}] &=&-\alpha e_{1}-\beta e_{2},  \notag \\
\lbrack e_{2},e_{3}] &=&\beta e_{1}+\alpha e_{2}+\alpha e_{3},\text{ \ \ }%
\alpha \neq 0.  \notag
\end{eqnarray}

In this case, $G=O(1,2)$ or $SL(2,R)$ if $\beta \neq 0$, while $G=E(1,1)$ if 
$\beta =0$.

- $\left( \mathfrak{g}_{2}\right) :$%
\begin{eqnarray}
\lbrack e_{1},e_{2}] &=&\gamma e_{2}-\beta e_{3},  \label{g2} \\
\lbrack e_{1},e_{3}] &=&-\beta e_{2}-\gamma e_{3},\text{ \ \ \ }\gamma \neq
0,  \notag \\
\lbrack e_{2},e_{3}] &=&\alpha e_{1}.  \notag
\end{eqnarray}

In this case, $G=O(1,2)$ or $SL(2,R)$ if $\alpha \neq 0$, while $G=E(1,1)$
if $\alpha =0$.

- $\left( \mathfrak{g}_{3}\right) :$%
\begin{eqnarray}
\lbrack e_{1},e_{2}] &=&-\gamma e_{3}  \label{g3} \\
\lbrack e_{1},e_{3}] &=&-\beta e_{2},  \notag \\
\lbrack e_{2},e_{3}] &=&\alpha e_{1}.  \notag
\end{eqnarray}

Table 1 lists all the Lie groups $G$ which admit a Lie algebra $\mathfrak{g}%
_{3}$, taking into account the different possibilities for $\alpha ,\beta $
and $\gamma $.%
\begin{eqnarray*}
&&%
\begin{tabular}{l||l||l||l}
$G$ & $\alpha $ & $\beta $ & $\gamma $ \\ \hline\hline
$O\left( 1,2\right) $ or $SL\left( 2,%
\mathbb{R}
\right) $ & $+$ & $+$ & $+$ \\ 
$O\left( 1,2\right) $ or $SL\left( 2,%
\mathbb{R}
\right) $ & $+$ & $-$ & $-$ \\ 
$SO\left( 3\right) $ or $SU\left( 2\right) $ & $+$ & $+$ & $-$ \\ 
$E\left( 2\right) $ & $+$ & $+$ & $0$ \\ 
$E\left( 2\right) $ & $+$ & $0$ & $-$ \\ 
$E\left( 1,1\right) $ & $+$ & $-$ & $0$ \\ 
$E\left( 1,1\right) $ & $+$ & $0$ & $+$ \\ 
$H_{3}$ & $+$ & $0$ & $0$ \\ 
$H_{3}$ & $0$ & $0$ & $-$ \\ 
$\mathbb{R} \oplus \mathbb{R} \oplus \mathbb{R}$ & $0$ & $0$ & $0$%
\end{tabular}
\\
&&\text{Table 1}
\end{eqnarray*}

- $\left( \mathfrak{g}_{4}\right) :$%
\begin{eqnarray}
\lbrack e_{1},e_{2}] &=&-e_{2}+(2\eta -\beta )e_{3},\text{ \ \ }\eta =\pm 1
\label{g4} \\
\lbrack e_{1},e_{3}] &=&-\beta e_{2}+e_{3},  \notag \\
\lbrack e_{2},e_{3}] &=&\alpha e_{1}.  \notag
\end{eqnarray}%
Table 2 describes all Lie groups $G$ admitting a Lie algebra $\mathfrak{g}%
_{4}.$%
\begin{eqnarray*}
&&%
\begin{tabular}{l||l||l}
$G$ & $\alpha $ & $\beta $ \\ \hline\hline
$O\left( 1,2\right) $ or $SL\left( 2,%
\mathbb{R}
\right) $ & $\neq 0$ & $\neq \eta $ \\ 
$E\left( 1,1\right) $ & $0$ & $\neq \eta $ \\ 
$E\left( 1,1\right) $ & $<0$ & $\eta $ \\ 
$E\left( 2\right) $ & $>0$ & $\eta $ \\ 
$H_{3}$ & $0$ & $\eta $%
\end{tabular}
\\
&&\text{Table 2}
\end{eqnarray*}
\end{theorem}

On the other hand, Cordero and Parker \cite{CP} classified three-dimensional
non-unimodular Lie groups equipped with a left-invariant Lorentzian metric,
obtaining a result corresponding to the one found by Milnor \cite{M76} in
the Riemannian case. In particular, they wrote down the possible forms of a
three-dimensional non-unimodular Lie algebra, determining their curvature
tensors and investigating the symmetry groups of the sectional curvature in
the different cases.

Let $G$ be a non-unimodular three-dimensional Lie group with a
left-invariant Lorentzian metric. Then the \textit{unimodular kernel} $%
\mathfrak{u}$ of the Lie algebra $\mathfrak{g}$\ of $G$ is defined by%
\begin{equation*}
\mathfrak{u}=\left\{ X\in \mathfrak{g}\text{ }/\text{ }\mathrm{tr}\text{ 
\textrm{ad}}\left( X\right) =0\right\} .
\end{equation*}%
Here ad $:\mathfrak{g}\rightarrow \mathrm{End}\left( \mathfrak{g}\right) $
is a homomorphism defined by \textrm{ad}$\left( X\right) Y=\left[ X,Y\right]
.$ One can see that \ is an ideal of $\mathfrak{g}$\ which contains the
ideal $\left[ \mathfrak{g,g}\right] .$

Following \cite{CP}, it was proven in \cite{C07} the following result:

\begin{theorem}
\label{non-unimodular} Let $(G,g)$ be a three-dimensional connected
non-unimodular Lie group, equipped with a left-invariant Lorentzian metric.
If $\ G$ is not symmetric, then, there exists a pseudo-orthonormal frame
field $\{e_{1},e_{2},e_{3}\}$, with $e_{3}$ timelike, such that the Lie
algebra of $G$ is one of the following:

- $\left( \mathfrak{g}_{5}\right) :$%
\begin{eqnarray}
\lbrack e_{1},e_{2}] &=&0,  \label{g5} \\
\lbrack e_{1},e_{3}] &=&\alpha e_{1}+\beta e_{2},  \notag \\
\lbrack e_{2},e_{3}] &=&\gamma e_{1}+\delta e_{2},\text{ \ \ }\alpha +\delta
\neq 0,\text{ }\alpha \gamma +\beta \delta =0.  \notag
\end{eqnarray}

- $\left( \mathfrak{g}_{6}\right) :$%
\begin{eqnarray}
\lbrack e_{1},e_{2}] &=&\alpha e_{2}+\beta e_{3},  \label{g6} \\
\lbrack e_{1},e_{3}] &=&\gamma e_{2}+\delta e_{3},\text{ \ \ \ }  \notag \\
\lbrack e_{2},e_{3}] &=&0,\text{ \ \ \ }\alpha +\delta \neq 0,\text{ }\alpha
\gamma -\beta \delta =0.  \notag
\end{eqnarray}

- $\left( \mathfrak{g}_{7}\right) :$%
\begin{eqnarray}
\lbrack e_{1},e_{2}] &=&-\alpha e_{1}-\beta e_{2}-\beta e_{3},  \label{g7} \\
\lbrack e_{1},e_{3}] &=&\alpha e_{1}+\beta e_{2}+\beta e_{3},\text{ \ \ \ } 
\notag \\
\lbrack e_{2},e_{3}] &=&\gamma e_{1}+\delta e_{2}+\delta e_{3},\text{ \ \ \ }%
\alpha +\delta \neq 0,\text{ }\alpha \gamma =0.  \notag
\end{eqnarray}
\end{theorem}

In the sequel, by $\left\{ G_{i}\right\} _{i=1,..,7}$ we shall denote the
connected, simply connected three-dimensional Lie group, equipped with a
left-invariant Lorentzian metric $g$ and having Lie algebra $\left\{ 
\mathfrak{g}_{i}\right\} _{i=1,..,7}.$ Let $\nabla $ be the Levi-Civita
connection of $G_{i}$ and $R$ its curvature tensor, taken with the sign
convention%
\begin{equation*}
R_{X,Y}=\nabla _{\left[ X,Y\right] }-\left[ \nabla _{X},\nabla _{Y}\right] .
\end{equation*}

The Ricci tensor of $\left( G_{i},g\right) $ is defined by%
\begin{equation*}
\varrho \left( X,Y\right) =\overset{3}{\underset{k=1}{\sum }}\varepsilon
_{k}g(R(X,e_{k})Y,e_{k}),
\end{equation*}%
where $\{e_{1},e_{2},e_{3}\}$ is a pseudo-orthonormal basis, with $e_{3}$
timelike and $\varepsilon _{1}=\varepsilon _{2}=1$ and $\varepsilon _{3}=-1,$
and the Ricci operator \textrm{Ric} is given by%
\begin{equation*}
\varrho \left( X,Y\right) =g\left( \mathrm{Ric}\left( X\right) ,Y\right)
\end{equation*}

\section{\protect\LARGE Algebraic Ricci Solitons of 3-dimensional unimodular Lie groups}

\setcounter{equation}{0}

\subsection{Algebraic Ricci Solitons of $G_{1}$}
We start by recalling the following result on the curvature tensor and the
Ricci operator of a three-dimensional Lorentzian Lie group $G_{1}.$

\begin{lemma}[\protect\cite{C07a}]
Let $\{e_{1},e_{2},e_{3}\}$ be the pseudo-orthonormal basis used in 
\eqref{g1}%
. Then%
\begin{eqnarray*}
\nabla _{e_{1}}e_{1} &=&-\alpha e_{2}-\alpha e_{3},\text{ \ \ \ }\nabla
_{e_{2}}e_{1}=\frac{\beta }{2}e_{3},\text{ \ \ \ \ \ \ \ \ \ \ \ \ \ }\nabla
_{e_{3}}e_{1}=\frac{\beta }{2}e_{2}, \\
\nabla _{e_{1}}e_{2} &=&\alpha e_{1}-\frac{\beta }{2}e_{3},\text{ \ \ \ \ }%
\nabla _{e_{2}}e_{2}=\alpha e_{3},\text{ \ \ \ \ \ \ \ \ \ \ \ \ \ \ }\nabla
_{e_{3}}e_{2}=-\frac{\beta }{2}e_{1}-\alpha e_{3}, \\
\nabla _{e_{1}}e_{3} &=&-\alpha e_{1}-\frac{\beta }{2}e_{2},\text{ \ \ }%
\nabla _{e_{2}}e_{3}=\frac{\beta }{2}e_{1}+\alpha e_{2},\text{ \ \ \ \ \ \ }%
\nabla _{e_{3}}e_{3}=-\alpha e_{2},
\end{eqnarray*}%
and the possible non-zero components of the curvature tensor are given by%
\begin{eqnarray*}
R_{1212} &=&-2\alpha ^{2}-\frac{\beta ^{2}}{4},\text{ \ \ \ \ }R_{1313}=%
\frac{\beta ^{2}}{4}-2\alpha ^{2},\text{ \ \ \ \ \ \ \ }R_{2323}=\frac{\beta
^{2}}{4}, \\
R_{1213} &=&2\alpha ^{2},\text{ \ \ \ \ \ \ \ \ \ \ \ \ \ }R_{1223}=-\alpha
\beta ,\text{ \ \ \ \ \ \ \ \ \ \ \ \ \ }R_{1323}=\alpha \beta .
\end{eqnarray*}%
Consequently, the Ricci operator is given by%
\begin{equation*}
\mathrm{Ric}=\left( 
\begin{array}{ccc}
-\frac{\beta ^{2}}{2} & -\alpha \beta & \alpha \beta \\ 
-\alpha \beta & -2\alpha ^{2}-\frac{\beta ^{2}}{2} & 2\alpha ^{2} \\ 
-\alpha \beta & -2\alpha ^{2} & 2\alpha ^{2}-\frac{\beta ^{2}}{2}%
\end{array}%
\right) .
\end{equation*}
\end{lemma}

Now, let $\mathrm{D}$ be an endomorphism of $\mathfrak{g}_{1}$ such that%
\begin{equation}
\mathrm{D}\left[ e_{i},e_{j}\right] =\left[ \mathrm{D}e_{i},e_{j}\right] +%
\left[ e_{i},\mathrm{D}e_{j}\right] \text{ for all }i,j  \label{D}
\end{equation}
where $\{e_{1},e_{2},e_{3}\}$ is the pseudo-orthonormal basis used in 
\eqref{g1}%
. Put 
\begin{equation*}
\mathrm{D}e_{l}=\lambda _{l}^{1}e_{1}+\lambda _{l}^{2}e_{2}+\lambda
_{l}^{3}e_{3}\text{ \ for all }l=1,2,3.
\end{equation*}

Starting from 
\eqref{g1}%
, we can write down 
\eqref{D}
and we get%
\begin{eqnarray}
&&%
\begin{array}{c}
\alpha \lambda _{2}^{2}-\alpha \lambda _{2}^{3}-\beta \lambda _{1}^{3}+\beta
\lambda _{3}^{1}=0,%
\end{array}
\label{sysD1} \\
&&%
\begin{array}{c}
\alpha \lambda _{1}^{3}+\beta \lambda _{2}^{3}+\alpha \lambda _{1}^{2}-\beta
\lambda _{3}^{2}=0,%
\end{array}
\notag \\
&&%
\begin{array}{c}
2\alpha \lambda _{1}^{3}+\beta \lambda _{1}^{1}+\beta \lambda _{2}^{2}-\beta
\lambda _{3}^{3}=0,%
\end{array}
\notag \\
&&%
\begin{array}{c}
-\alpha \lambda _{3}^{3}+\beta \lambda _{1}^{2}+\alpha \lambda
_{3}^{2}+\beta \lambda _{2}^{1}=0,%
\end{array}
\notag \\
&&%
\begin{array}{c}
2\alpha \lambda _{1}^{2}-\beta \lambda _{1}^{1}-\beta \lambda _{3}^{3}+\beta
\lambda _{2}^{2}=0,%
\end{array}
\notag \\
&&%
\begin{array}{c}
-2\alpha \lambda _{2}^{1}-2\alpha \lambda _{3}^{1}+\beta \lambda
_{2}^{2}+\beta \lambda _{3}^{3}-\beta \lambda _{1}^{1}=0.%
\end{array}
\notag
\end{eqnarray}

The following result can now verified by a straightforward calculation.

\begin{lemma}
The solutions of 
\eqref{sysD1}
are the following:

\begin{itemize}
\item If $\beta \neq 0:$%
\begin{equation*}
\mathrm{D=}\left( 
\begin{array}{ccc}
\lambda _{1}^{1} & \lambda _{2}^{1} & -\frac{\beta }{\alpha }\lambda
_{1}^{1}-\lambda _{2}^{1} \\ 
-\frac{\beta }{\alpha }\lambda _{2}^{2} & \lambda _{2}^{2} & -\lambda
_{1}^{1}+\left( \frac{\beta ^{2}}{\alpha ^{2}}-1\right) \lambda _{2}^{2}-%
\frac{\beta }{\alpha }\lambda _{2}^{1} \\ 
-\frac{\beta }{\alpha }\left( \lambda _{1}^{1}+\lambda _{2}^{2}\right) & 
\left( \frac{\beta ^{2}}{\alpha ^{2}}+1\right) \lambda _{2}^{2}-\frac{\beta 
}{\alpha }\lambda _{2}^{1} & -\lambda _{1}^{1}-\lambda _{2}^{2}%
\end{array}%
\right) .
\end{equation*}

\item If $\beta =0:$%
\begin{equation*}
\mathrm{D=}\left( 
\begin{array}{ccc}
\lambda _{1}^{1} & \lambda _{2}^{1} & -\lambda _{2}^{1} \\ 
0 & \lambda _{2}^{2} & \lambda _{3}^{3} \\ 
0 & \lambda _{2}^{2} & \lambda _{3}^{3}%
\end{array}%
\right) .
\end{equation*}
\end{itemize}
\end{lemma}

We now prove the following.

\begin{theorem}
\label{ASG1}Consider the three-dimensional Lorentzian Lie group $G_{1}$.
Then, $G_{1}$ is an algebraic Ricci Soliton Lie group if and only if $\beta
=0.$ In particular%
\begin{equation*}
\mathrm{D=Ric=}\left( 
\begin{array}{ccc}
0 & 0 & 0 \\ 
0 & -2\alpha ^{2} & 2\alpha ^{2} \\ 
0 & -2\alpha ^{2} & 2\alpha ^{2}%
\end{array}%
\right) \text{ and }c=0.
\end{equation*}
\end{theorem}

\textbf{Proof.} $G_{1}$ is an algebraic Ricci Soliton Lie group if and only
if 
\eqref{Soliton}
holds. Hence, we have to consider two cases:

-If $\beta \neq 0,$ in this case 
\eqref{Soliton}
is equivalent to the following system 
\begin{eqnarray}
&&%
\begin{array}{c}
c+\lambda _{1}^{1}=-\frac{\beta ^{2}}{2},\text{ \ }%
\end{array}
\label{Sol1} \\
&&%
\begin{array}{c}
c+\lambda _{2}^{2}=-2\alpha ^{2}-\frac{\beta ^{2}}{2},%
\end{array}
\notag \\
&&%
\begin{array}{c}
\text{\ }c-\lambda _{1}^{1}-\lambda _{2}^{2}=2\alpha ^{2}-\frac{\beta ^{2}}{2%
},%
\end{array}
\notag \\
&&%
\begin{array}{c}
\lambda _{2}^{2}=\alpha ^{2},\text{ }%
\end{array}
\notag \\
&&%
\begin{array}{c}
\left( \lambda _{1}^{1}+\lambda _{2}^{2}\right) =\alpha ^{2},%
\end{array}
\notag \\
&&%
\begin{array}{c}
\lambda _{2}^{1}=-\alpha \beta ,%
\end{array}
\notag \\
&&%
\begin{array}{c}
-\frac{\beta }{\alpha }\lambda _{1}^{1}-\lambda _{2}^{1}=\alpha \beta ,%
\end{array}
\notag \\
&&%
\begin{array}{c}
\left( \frac{\beta ^{2}}{\alpha ^{2}}+1\right) \lambda _{2}^{2}-\frac{\beta 
}{\alpha }\lambda _{2}^{1}=-2\alpha ^{2},\text{ \ }%
\end{array}%
\text{\ }  \notag \\
&&%
\begin{array}{c}
-\lambda _{1}^{1}+\left( \frac{\beta ^{2}}{\alpha ^{2}}-1\right) \lambda
_{2}^{2}-\frac{\beta }{\alpha }\lambda _{2}^{1}=2\alpha ^{2},%
\end{array}%
\text{\ \ \ \ \ \ \ \ }  \notag
\end{eqnarray}%
for some real constant $c.$ We conclude that $\lambda _{1}^{1}=0$ and $c=-%
\frac{\beta ^{2}}{2}$. Hence one obtain, from the third equation in 
\eqref{Sol1}%
, that $\lambda _{2}^{2}=-2\alpha ^{2}=\alpha ^{2}.$ Therefore 
\eqref{Sol1}
does not occur for $\alpha \neq 0.$

-If $\beta =0,$ we find that 
\eqref{Soliton}
is satisfied if and only if%
\begin{equation*}
\lambda _{3}^{3}=-\lambda _{2}^{2}=2\alpha ^{2}\text{ \ and }\lambda
_{1}^{1}=\lambda _{2}^{1}=c=0. \hspace{3cm} \square
\end{equation*}

\begin{remark} 
Let $D$ be the derivation given in the above Theorem. 
Then D is an inner derivation, that is, $D=2\alpha (\mathrm{ad}e_{2}+\mathrm{ad}e_{3}).$
\end{remark}

\begin{remark}
Let $\beta =0$ then $G_{1}$ is the solvable Lie group $E(1,1)$. So we can
construct a pseudo-orthonormal basis $\left\{ e_{1},e_{2},e_{3}\right\} $,
by setting%
\begin{eqnarray*}
e_{1} &=&\dfrac{1}{2}e^{-2\alpha x_{1}}\Big(\dfrac{\partial }{\partial x_{2}}%
-\dfrac{\partial }{\partial x_{3}}\Big), \\
e_{2} &=&\dfrac{1}{2}\dfrac{\partial }{\partial x_{1}}+\dfrac{1}{4}%
e^{2\alpha x_{1}}\Big(\dfrac{\partial }{\partial x_{2}}+\dfrac{\partial }{%
\partial x_{3}}\Big), \\
e_{3} &=&-\dfrac{1}{2}\dfrac{\partial }{\partial x_{1}}+\dfrac{1}{4}%
e^{2\alpha x_{1}}\Big(\dfrac{\partial }{\partial x_{2}}+\dfrac{\partial }{%
\partial x_{3}}\Big).
\end{eqnarray*}

Using Theorem \ref{ASG1} we obtain that%
\begin{equation*}
\mathrm{\exp }\left( \frac{t}{2}\mathrm{D}\right) =\left( 
\begin{array}{ccc}
1 & 0 & 0 \\ 
0 & 1-t\alpha ^{2} & t\alpha ^{2} \\ 
0 & -t\alpha ^{2} & 1+t\alpha ^{2}%
\end{array}%
\right) .
\end{equation*}%
Next, let $\varphi =(\varphi ^{1},\varphi ^{2},\varphi ^{3})$ such that 
\begin{equation}
\mathrm{\exp }\left( \frac{t}{2}\mathrm{D}\right) =d\varphi _{t}|_{e}.
\label{exp1}
\end{equation}%
Standard computations lead to conclude that 
\eqref{exp1}
is satisfied if and only if the following system holds%
\begin{eqnarray*}
&&%
\begin{array}{c}
\varphi _{1}^{2}+\dfrac{e^{2\alpha x_{1}}}{2}(\varphi _{2}^{2}+\varphi
_{3}^{2})=\dfrac{e^{2\alpha x_{1}}}{2}-t\alpha ^{2}e^{2\alpha x_{1}},%
\end{array}
\\
&&%
\begin{array}{c}
\varphi _{1}^{3}+\dfrac{e^{2\alpha x_{1}}}{2}(\varphi _{2}^{3}+\varphi
_{3}^{3})=\dfrac{e^{2\alpha x_{1}}}{2}-t\alpha ^{2}e^{2\alpha x_{1}},%
\end{array}
\\
&&%
\begin{array}{c}
-\varphi _{1}^{2}+\dfrac{e^{2\alpha x_{1}}}{2}(\varphi _{2}^{2}+\varphi
_{3}^{2})=\dfrac{e^{2\alpha x_{1}}}{2}+t\alpha ^{2}e^{2\alpha x_{1}},%
\end{array}
\\
&&%
\begin{array}{c}
-\varphi _{1}^{3}+\dfrac{e^{2\alpha x_{1}}}{2}(\varphi _{2}^{3}+\varphi
_{3}^{3})=\dfrac{e^{2\alpha x_{1}}}{2}+t\alpha ^{2}e^{2\alpha x_{1}},%
\end{array}
\\
&&%
\begin{array}{c}
\varphi _{1}^{1}+\dfrac{e^{2\alpha x_{1}}}{2}(\varphi _{2}^{1}+\varphi
_{3}^{1})=1,%
\end{array}
\\
&&%
\begin{array}{c}
-\varphi _{1}^{1}+\dfrac{e^{2\alpha x_{1}}}{2}(\varphi _{2}^{1}+\varphi
_{3}^{1})=-1,%
\end{array}
\\
&&%
\begin{array}{c}
\varphi _{2}^{1}-\varphi _{3}^{1}=0,%
\end{array}
\\
&&%
\begin{array}{c}
\varphi _{2}^{2}-\varphi _{3}^{2}=1,%
\end{array}
\\
&&%
\begin{array}{c}
\varphi _{2}^{3}-\varphi _{3}^{3}=-1,%
\end{array}%
\end{eqnarray*}%
where $\varphi _{i}^{j}=\dfrac{\partial \varphi ^{j}}{\partial x_{i}}$.
Hence, we deduce that 
\begin{equation*}
\varphi =\left( x_{1},x_{2}-\dfrac{t}{2}\alpha e^{2\alpha x_{1}},x_{3}-%
\dfrac{t}{2}\alpha e^{2\alpha x_{1}}\right) .
\end{equation*}%
The vector field, for which 
\eqref{Ricci soliton}
holds, is then given by 
\begin{equation*}
X_{D}=\dfrac{d\varphi _{t}}{dt}\Big|_{t=0}(p)=\left( 0,-\dfrac{\alpha }{2}%
e^{2\alpha x_{1}},-\dfrac{\alpha }{2}e^{2\alpha x_{1}}\right) =-\alpha
(e_{2}+e_{3}).
\end{equation*}%
$X_{D}$ coincide with the Ricci soliton vector obtained by M. Brozos-Vazquez
et. al (see Theorem 1 case (iv.2) of \cite{BCGG09}).
\end{remark}

\begin{remark}
\label{maintheorem2a} Now, let $\beta \neq 0$ then $G_{1}$ is the Lie group $%
SL(2,%
\mathbb{R}
).$ Put $Y=\sum Y^{i}e_{i}$ where $\{e_{1},e_{2},e_{3}\}$ denotes the basis
used in 
\eqref{g1}
and let $g$ be the left-invariant Lorentzian metric for which $%
\{e_{1},e_{2},e_{3}\}$ is the pseudo-orthonormal for $G_{1}$. Therefore%
\begin{equation*}
L_{Y}g=\left( 
\begin{array}{ccc}
2\alpha \left( Y_{2}-Y_{3}\right) & -\alpha Y_{1} & \alpha Y_{1} \\ 
-\alpha Y_{1} & 2\alpha Y_{3} & -\alpha \left( Y_{2}+Y_{3}\right) \\ 
\alpha Y_{1} & -\alpha \left( Y_{2}+Y_{3}\right) & 2\alpha Y_{2}%
\end{array}%
\right) ,\alpha \neq 0,
\end{equation*}%
and hence, we have a Ricci soliton equation 
\eqref{Ricci soliton}
is satisfied if and only if%
\begin{eqnarray*}
&&%
\begin{array}{c}
2\alpha \left( Y_{2}-Y_{3}\right) +c=-\frac{\beta ^{2}}{2},%
\end{array}
\\
&&%
\begin{array}{c}
2\alpha Y_{3}+c=-2\alpha ^{2}-\frac{\beta ^{2}}{2},%
\end{array}
\\
&&%
\begin{array}{c}
2\alpha Y_{2}-c=-2\alpha ^{2}+\frac{\beta ^{2}}{2},%
\end{array}
\\
&&%
\begin{array}{c}
\alpha \left( Y_{2}+Y_{3}\right) =-2\alpha ^{2},%
\end{array}
\\
&&%
\begin{array}{c}
Y_{1}=\beta .%
\end{array}%
\end{eqnarray*}%
Therefore, one easily gets that $Y_{1}=\beta ,$ $Y_{2}=Y_{3}=-\alpha $ and
that $c=-\frac{\beta ^{2}}{2};$ thus, there exist a left-invariant Ricci
soliton vector field, given by%
\begin{equation*}
Y=\beta e_{1}-\alpha \left( e_{2}+e_{3}\right) \text{.}
\end{equation*}%
Note that, $G_{1}=SL\left( 2,%
\mathbb{R}
\right) $ is a Ricci soliton however there do no exist any algebraic Ricci soliton in $G_{1}$ with $\beta \neq 0.$ 
\end{remark}

\subsection{Algebraic Ricci Solitons of $G_{2}$}
Next we consider $G_{2}.$ The following lemma was proven in \cite{C07a}, but
we correct here a misprint in some of the Levi-Civita, curvature and the
Ricci operator components of a three-dimensional Lorentzian Lie group $%
G_{2}. $

\begin{lemma}
Let $\{e_{1},e_{2},e_{3}\}$ be the pseudo-orthonormal basis used in 
\eqref{g2}%
. Then%
\begin{eqnarray*}
\nabla _{e_{1}}e_{1} &=&0,\text{ \ \ \ \ \ \ \ \ \ \ \ \ \ \ \ \ \ \ \ \ }%
\nabla _{e_{2}}e_{1}=-\gamma e_{2}+\frac{\alpha }{2}e_{3},\text{ \ \ \ \ \ \
\ \ \ }\nabla _{e_{3}}e_{1}=\frac{\alpha }{2}e_{2}+\gamma e_{3}, \\
\nabla _{e_{1}}e_{2} &=&\left( \frac{\alpha }{2}-\beta \right) e_{3},\text{
\ \ \ \ \ \ \ }\nabla _{e_{2}}e_{2}=\gamma e_{1},\text{ \ \ \ \ \ \ \ \ \ \
\ \ \ \ \ \ \ \ \ \ }\nabla _{e_{3}}e_{2}=-\frac{\alpha }{2}e_{1}, \\
\nabla _{e_{1}}e_{3} &=&\left( \frac{\alpha }{2}-\beta \right) e_{2},\text{
\ \ \ \ \ \ \ }\nabla _{e_{2}}e_{3}=\frac{\alpha }{2}e_{1},\text{ \ \ \ \ \
\ \ \ \ \ \ \ \ \ \ \ \ \ \ }\nabla _{e_{3}}e_{3}=\gamma e_{1},
\end{eqnarray*}%
and the possible non-vanishing components of the curvature tensor are given
by%
\begin{equation*}
R_{1212}=-\gamma ^{2}-\frac{\alpha ^{2}}{4},\text{ \ }R_{1213}=\gamma \left(
2\beta -\alpha \right) ,\text{ \ }R_{1313}=\gamma ^{2}+\frac{\alpha ^{2}}{4},%
\text{\ \ \ }R_{2323}=-\gamma ^{2}-\frac{3\alpha ^{2}}{4}+\alpha \beta .
\end{equation*}%
The Ricci operator is given by%
\begin{equation*}
\mathrm{Ric}=\left( 
\begin{array}{ccc}
-\frac{\alpha ^{2}}{2}-2\gamma ^{2} & 0 & 0 \\ 
0 & \alpha \left( \frac{\alpha }{2}-\beta \right) & \gamma \left( 2\beta
-\alpha \right) \\ 
0 & \gamma \left( \alpha -2\beta \right) & \alpha \left( \frac{\alpha }{2}%
-\beta \right)%
\end{array}%
\right) .
\end{equation*}
\end{lemma}

Next, let $\mathrm{D}e_{l}=\lambda _{l}^{1}e_{1}+\lambda
_{l}^{2}e_{2}+\lambda _{l}^{3}e_{3}$ be an endomorphism of $\mathfrak{g}_{2}$
where $\{e_{1},e_{2},e_{3}\}$ is the pseudo-orthonormal basis used in 
\eqref{g2}%
. Then, $\mathrm{D}\in \mathrm{Der}\left( \mathfrak{g}_{2}\right) $ if and
only if%
\begin{eqnarray}
&&%
\begin{array}{c}
\gamma \lambda _{2}^{1}-\beta \lambda _{3}^{1}+\alpha \lambda _{1}^{3}=0,%
\end{array}
\label{sysD2} \\
&&%
\begin{array}{c}
\gamma \lambda _{1}^{1}-\beta \lambda _{2}^{3}+\beta \lambda _{3}^{2}=0,%
\end{array}
\notag \\
&&%
\begin{array}{c}
\beta \lambda _{1}^{1}+\beta \lambda _{2}^{2}+2\gamma \lambda _{2}^{3}-\beta
\lambda _{3}^{3}=0,%
\end{array}
\notag \\
&&%
\begin{array}{c}
\beta \lambda _{2}^{1}+\gamma \lambda _{3}^{1}+\alpha \lambda _{1}^{2}=0,%
\end{array}
\notag \\
&&%
\begin{array}{c}
\beta \lambda _{1}^{1}+\beta \lambda _{3}^{3}-2\gamma \lambda _{3}^{2}-\beta
\lambda _{2}^{2}=0,%
\end{array}
\notag \\
&&%
\begin{array}{c}
\alpha \left( \lambda _{2}^{2}+\lambda _{3}^{3}-\lambda _{1}^{1}\right) =0.%
\end{array}
\notag
\end{eqnarray}

Therefore, the following result holds:

\begin{lemma}
The solutions of 
\eqref{sysD2}
are the following:

\begin{itemize}
\item If $\alpha \neq 0:$%
\begin{equation*}
\mathrm{D=}\left( 
\begin{array}{ccc}
0 & \lambda _{2}^{1} & \lambda _{3}^{1} \\ 
-\frac{1}{\alpha }\left( \gamma \lambda _{3}^{1}+\beta \lambda
_{2}^{1}\right) & \lambda _{2}^{2} & -\frac{\beta }{\gamma }\lambda _{2}^{2}
\\ 
\frac{1}{\alpha }\left( \beta \lambda _{3}^{1}-\gamma \lambda _{2}^{1}\right)
& -\frac{\beta }{\gamma }\lambda _{2}^{2} & -\lambda _{2}^{2}%
\end{array}%
\right) .
\end{equation*}

\item If $\alpha =0:$%
\begin{equation*}
\mathrm{D=}\left( 
\begin{array}{ccc}
0 & 0 & 0 \\ 
\lambda _{1}^{2} & \lambda _{2}^{2} & \frac{\beta }{2\gamma }\left( -\lambda
_{2}^{2}+\lambda _{3}^{3}\right) \\ 
\lambda _{1}^{3} & \frac{\beta }{2\gamma }\left( -\lambda _{2}^{2}+\lambda
_{3}^{3}\right) & \lambda _{3}^{3}%
\end{array}%
\right) .
\end{equation*}
\end{itemize}
\end{lemma}

Using the above lemma, we now prove the following.

\begin{theorem}
\label{ASG2}Consider the three-dimensional Lorentzian Lie group $G_{2}$.
Then, $G_{2}$ is an algebraic Ricci Soliton Lie group if and only if $\alpha
=\beta =0.$ In particular%
\begin{equation*}
\mathrm{D=}\left( 
\begin{array}{ccc}
0 & 0 & 0 \\ 
0 & 2\gamma ^{2} & 0 \\ 
0 & 0 & 2\gamma ^{2}%
\end{array}%
\right) \text{ and }c=-2\gamma ^{2}.
\end{equation*}
\end{theorem}

\textbf{Proof. }We need to distinguish two cases:

-If $\alpha \neq 0,$ in this case the algebraic Ricci Soliton condition 
\eqref{Soliton}
on $G_{2}$ gives rise to the following system:%
\begin{eqnarray}
&&%
\begin{array}{c}
c=-\left( \frac{\alpha ^{2}}{2}+2\gamma ^{2}\right) ,\text{ \ }%
\end{array}
\label{Sol2} \\
&&%
\begin{array}{c}
c+\lambda _{2}^{2}=c-\lambda _{2}^{2}=\frac{\alpha ^{2}}{2}-\alpha \beta ,%
\end{array}
\notag \\
&&%
\begin{array}{c}
-\beta \lambda _{2}^{2}=\gamma ^{2}\left( \alpha -2\beta \right) =\gamma
^{2}\left( 2\beta -\alpha \right) ,%
\end{array}
\notag \\
&&%
\begin{array}{c}
\lambda _{2}^{1}=\lambda _{3}^{1}=0.%
\end{array}
\notag
\end{eqnarray}

By the second equation in 
\eqref{Sol2}%
, we get%
\begin{equation*}
c=\frac{\alpha ^{2}}{2}-\alpha \beta \text{ \ and \ }\lambda _{2}^{2}=0.
\end{equation*}%
Hence, using the third equation in 
\eqref{Sol2}
we obtain, since $\gamma \neq 0,$ that%
\begin{equation*}
\alpha =2\beta .
\end{equation*}%
Therefore, it follows from the third equation of 
\eqref{Sol2}
that $c=0,$ and the first equation in 
\eqref{Sol2}
becomes $\frac{\alpha ^{2}}{2}+2\gamma ^{2}=0$ which is a contradiction,
since $\gamma \neq 0.$

-If $\alpha =0,$ we find that 
\eqref{Soliton}
is satisfied if and only if:%
\begin{eqnarray}
&&%
\begin{array}{c}
c=-2\gamma ^{2},\text{ \ }%
\end{array}
\label{Sol2a} \\
&&%
\begin{array}{c}
c+\lambda _{2}^{2}=c+\lambda _{3}^{3}=0,%
\end{array}
\notag \\
&&%
\begin{array}{c}
\beta \left( -\lambda _{2}^{2}+\lambda _{3}^{3}\right) =4\beta \gamma
^{2}=-4\beta \gamma ^{2},%
\end{array}
\notag \\
&&%
\begin{array}{c}
\lambda _{1}^{2}=\lambda _{1}^{3}=0.%
\end{array}
\notag
\end{eqnarray}%
It is easy to show that, whenever $\beta \neq 0,$ 
\eqref{Soliton}
does not occur (since $\gamma \neq 0$). On the other hand, if $\beta =0$
then 
\begin{equation*}
\lambda _{2}^{2}=\lambda _{3}^{3}=-c=2\gamma ^{2}\text{ \ and }\lambda
_{1}^{2}=\lambda _{1}^{3}=0,
\end{equation*}%
is the only solution of 
\eqref{Soliton}%
. $\square$

\begin{remark} 
Let $D$ be the derivation given in the above Theorem. Routine calculations
show that D is an outer derivation, that is, $D\neq \mathrm{ad}X$, for any $X\in \mathfrak{g}_{2}$.
\end{remark}

\begin{remark}
Let $\alpha =\beta =0$, then $G_{2}$ is the solvable Lie group $E(1,1)$. We
construct a pseudo-orthonormal basis $\left\{ e_{1},e_{2},e_{3}\right\} $,
by setting%
\begin{eqnarray*}
e_{1} &=&\dfrac{\partial }{\partial x_{1}}, \\
e_{2} &=&\dfrac{1}{2}e^{\gamma x_{1}}\Big(\dfrac{\partial }{\partial x_{2}}+%
\dfrac{\partial }{\partial x_{3}}\Big), \\
e_{3} &=&\dfrac{1}{2}e^{-\gamma x_{1}}\Big(\dfrac{\partial }{\partial x_{2}}-%
\dfrac{\partial }{\partial x_{3}}\Big).
\end{eqnarray*}%
Using Theorem \ref{ASG2} we obtain that%
\begin{equation*}
\mathrm{\exp }\left( \frac{t}{2}\mathrm{D}\right) =\left( 
\begin{array}{ccc}
1 & 0 & 0 \\ 
0 & \exp \left( \gamma ^{2}t\right) & 0 \\ 
0 & 0 & \exp \left( \gamma ^{2}t\right)%
\end{array}%
\right) .
\end{equation*}%
Next, let $\psi =(\psi ^{1},\psi ^{2},\psi ^{3})$ such that 
\begin{equation}
\mathrm{\exp }\left( \frac{t}{2}\mathrm{D}\right) =d\psi _{t}|_{e}.
\label{exp2}
\end{equation}%
Standard computations show that 
\eqref{exp2}
is satisfied if and only if the following system holds:%
\begin{eqnarray*}
&&%
\begin{array}{c}
\psi _{2}^{2}+\psi _{3}^{2}=e^{\gamma ^{2}t},%
\end{array}
\\
&&%
\begin{array}{c}
\psi _{2}^{3}+\psi _{3}^{3}=e^{\gamma ^{2}t},%
\end{array}
\\
&&%
\begin{array}{c}
\psi _{2}^{2}-\psi _{3}^{2}=e^{\gamma ^{2}t},%
\end{array}
\\
&&%
\begin{array}{c}
\psi _{2}^{3}-\psi _{3}^{3}=-e^{\gamma ^{2}t},%
\end{array}
\\
&&%
\begin{array}{c}
\psi _{1}^{1}=1,\text{ }\psi _{1}^{2}=\psi _{1}^{3}=\psi _{2}^{1}=\psi
_{3}^{1}=0.%
\end{array}%
\end{eqnarray*}%
where $\psi _{i}^{j}=\dfrac{\partial \psi ^{j}}{\partial x_{i}}$. So, we
obtain 
\begin{equation*}
\psi =\left( x_{1},e^{\gamma ^{2}t}x_{2},e^{\gamma ^{2}t}x_{3}\right) .
\end{equation*}%
Therefore, the vector field, for which 
\eqref{Ricci soliton}
holds, is given by 
\begin{eqnarray*}
X_{D} &=&\dfrac{d\psi _{t}}{dt}\Big|_{t=0}(p) \\
&=&\left( 0,\gamma ^{2}e^{\gamma ^{2}t}x_{2},\gamma ^{2}e^{\gamma
^{2}t}x_{3}\right) \\
&=&\gamma ^{2}e^{-\gamma x_{1}}(x_{2}+x_{3})e_{2}+\gamma ^{2}e^{\gamma
x_{1}}(x_{2}-x_{3})e_{3}.
\end{eqnarray*}%
Notice that $X_{D}$ is not a left-invariant vector field, this means that $%
X_{D}$ is a new Ricci soliton vector for $G_{2}$ ($\alpha =\beta =0$).
\end{remark}


\subsection{Algebraic Ricci Solitons of $G_{3}$}
Next consider now the Lorentzian Lie group $G_{3},$ we start from the
following result of \cite{C07a}.

\begin{lemma}
Let $\{e_{1},e_{2},e_{3}\}$ be the pseudo-orthonormal basis used in 
\eqref{g3}%
. Then%
\begin{eqnarray*}
\nabla _{e_{1}}e_{1} &=&0,\text{ \ \ \ \ \ \ \ \ }\nabla
_{e_{2}}e_{1}=a_{2}e_{3},\text{ \ \ \ \ }\nabla _{e_{3}}e_{1}=a_{3}e_{2}, \\
\nabla _{e_{1}}e_{2} &=&a_{1}e_{3},\text{ \ \ \ \ }\nabla _{e_{2}}e_{2}=0,%
\text{\ \ \ \ \ \ \ \ \ }\nabla _{e_{3}}e_{2}=-a_{3}e_{1}, \\
\nabla _{e_{1}}e_{3} &=&a_{1}e_{2},\text{ \ \ \ \ }\nabla
_{e_{2}}e_{3}=a_{2}e_{1},\text{ \ \ \ \ }\nabla _{e_{3}}e_{3}=0,
\end{eqnarray*}%
where we put%
\begin{equation*}
a_{1}=\frac{1}{2}\left( \alpha -\beta -\gamma \right) ,\text{ \ \ \ }a_{2}=%
\frac{1}{2}\left( \alpha -\beta +\gamma \right) ,\text{ \ \ \ }a_{3}=\frac{1%
}{2}\left( \alpha +\beta -\gamma \right) .
\end{equation*}%
The not always null components of the curvature tensor are given by%
\begin{equation*}
R_{1212}=-\left( a_{1}a_{2}+\gamma a_{3}\right) ,\text{ \ }R_{1313}=\left(
a_{1}a_{3}+\beta a_{2}\right) ,\text{ \ \ }R_{2323}=-\left(
a_{2}a_{3}+\alpha a_{3}\right) .
\end{equation*}%
The Ricci operator is given by%
\begin{equation*}
\mathrm{Ric}=\left( 
\begin{array}{ccc}
-a_{1}a_{2}-a_{1}a_{3}-\beta a_{2}-\gamma a_{3} & 0 & 0 \\ 
0 & a_{2}a_{3}-a_{1}a_{2}+\alpha a_{1}-\gamma a_{3} & 0 \\ 
0 & 0 & -a_{1}a_{3}+a_{2}a_{3}+\alpha a_{1}-\beta a_{3}%
\end{array}%
\right) .
\end{equation*}
\end{lemma}

Next, put $\mathrm{D}e_{l}=\lambda _{l}^{1}e_{1}+\lambda
_{l}^{2}e_{2}+\lambda _{l}^{3}e_{3}$ where $\{e_{1},e_{2},e_{3}\}$ is the
pseudo-orthonormal basis used in 
\eqref{g3}%
. Then, $\mathrm{D}\in \mathrm{Der}\left( \mathfrak{g}_{3}\right) $ if and
only if%
\begin{eqnarray}
&&%
\begin{array}{c}
\beta \lambda _{2}^{3}-\gamma \lambda _{3}^{2}=0,%
\end{array}
\label{sysD3} \\
&&%
\begin{array}{c}
\gamma \left( \lambda _{1}^{1}+\lambda _{2}^{2}-\lambda _{3}^{3}\right) =0,%
\end{array}
\notag \\
&&%
\begin{array}{c}
\alpha \lambda _{1}^{2}+\beta \lambda _{2}^{1}=0,%
\end{array}
\notag \\
&&%
\begin{array}{c}
\beta \left( \lambda _{1}^{1}-\lambda _{2}^{2}+\lambda _{3}^{3}\right) =0,%
\end{array}
\notag \\
&&%
\begin{array}{c}
\alpha \left( -\lambda _{1}^{1}+\lambda _{2}^{2}+\lambda _{3}^{3}\right) =0,%
\end{array}
\notag \\
&&%
\begin{array}{c}
\gamma \lambda _{3}^{1}-\alpha \lambda _{1}^{3}=0.%
\end{array}
\notag
\end{eqnarray}

Routine, but long calculations lead to the following:

\begin{lemma}
The solutions of 
\eqref{sysD3}
are the following:

\begin{itemize}
\item If $\beta =\gamma =0\neq \alpha :$%
\begin{equation*}
\mathrm{D=}\left( 
\begin{array}{ccc}
\lambda _{2}^{2}+\lambda _{3}^{3} & \lambda _{2}^{1} & \lambda _{3}^{1} \\ 
0 & \lambda _{2}^{2} & \lambda _{3}^{2} \\ 
0 & \lambda _{2}^{3} & \lambda _{3}^{3}%
\end{array}%
\right) .
\end{equation*}

\item If $\alpha =\beta =0\neq \gamma :$%
\begin{equation*}
\mathrm{D=}\left( 
\begin{array}{ccc}
\lambda _{3}^{3}-\lambda _{2}^{2} & \lambda _{2}^{1} & 0 \\ 
\lambda _{1}^{2} & \lambda _{2}^{2} & 0 \\ 
\lambda _{1}^{3} & \lambda _{2}^{3} & \lambda _{3}^{3}%
\end{array}%
\right) .
\end{equation*}

\item If $\alpha \beta \gamma \neq 0:$ 
\begin{equation*}
\mathrm{D=}\left( 
\begin{array}{ccc}
0 & -\frac{\alpha }{\beta }\lambda _{1}^{1} & \frac{\alpha }{\gamma }\lambda
_{1}^{3} \\ 
\lambda _{1}^{2} & 0 & \frac{\beta }{\gamma }\lambda _{2}^{3} \\ 
\lambda _{1}^{3} & \lambda _{2}^{3} & 0%
\end{array}%
\right) .
\end{equation*}

\item If $\alpha \beta \neq 0=\gamma :$%
\begin{equation*}
\mathrm{D=}\left( 
\begin{array}{ccc}
\lambda _{1}^{1} & -\frac{\alpha }{\beta }\lambda _{1}^{2} & \lambda _{3}^{1}
\\ 
\lambda _{1}^{2} & \lambda _{1}^{1} & \lambda _{3}^{2} \\ 
0 & 0 & 0%
\end{array}%
\right) .
\end{equation*}

\item If $\alpha \gamma \neq 0=\beta :$%
\begin{equation*}
\mathrm{D=}\left( 
\begin{array}{ccc}
\lambda _{1}^{1} & \lambda _{2}^{1} & \frac{\alpha }{\gamma }\lambda _{1}^{3}
\\ 
0 & 0 & 0 \\ 
\lambda _{1}^{3} & \lambda _{2}^{3} & \lambda _{1}^{1}%
\end{array}%
\right) .
\end{equation*}

\item If $\alpha =\beta =\gamma =0$, in this case $\mathrm{D}e_{l}=\lambda
_{l}^{1}e_{1}+\lambda _{l}^{2}e_{2}+\lambda _{l}^{3}e_{3}$.
\end{itemize}
\end{lemma}

We can now prove the following:

\begin{theorem}
Consider the three-dimensional Lorentzian Lie group $G_{3}$. Then, $G_{3}$
is an algebraic Ricci Soliton Lie group if and only if one of the following
statements holds true:

\begin{itemize}
\item $\alpha >0$ and $\beta =\gamma =0.$ In particular%
\begin{equation*}
\mathrm{D=}\left( 
\begin{array}{ccc}
-2\alpha ^{2} & 0 & 0 \\ 
0 & -\alpha ^{2} & 0 \\ 
0 & 0 & -\alpha ^{2}%
\end{array}%
\right) \text{ and }c=\frac{3\alpha ^{2}}{2}.
\end{equation*}

\item $\gamma <0$ and $\alpha =\beta =0.$ In this case we have%
\begin{equation*}
\mathrm{D=}\left( 
\begin{array}{ccc}
-\gamma ^{2} & 0 & 0 \\ 
0 & -\gamma ^{2} & 0 \\ 
0 & 0 & -2\gamma ^{2}%
\end{array}%
\right) \text{ and }c=\frac{3\gamma ^{2}}{2}.
\end{equation*}

\item $\alpha =\beta =\gamma \neq 0$. In this case we have $\mathrm{D}=0$
and $c=-\frac{\alpha ^{2}}{2}.$

\item $\alpha =\beta >0$ and $\gamma =0$. In this case we have $\mathrm{D}=0$
and $c=0.$

\item $\beta =-\alpha <0$ and $\gamma =0$. In particular%
\begin{equation*}
\mathrm{D=}\left( 
\begin{array}{ccc}
-2\alpha ^{2} & 0 & 0 \\ 
0 & -2\alpha ^{2} & 0 \\ 
0 & 0 & 0%
\end{array}%
\right) \text{ and }c=2\alpha ^{2}.
\end{equation*}

\item $\gamma =-\alpha <0$ and $\beta =0$. In this case we have%
\begin{equation*}
\mathrm{D=}\left( 
\begin{array}{ccc}
-2\alpha ^{2} & 0 & 0 \\ 
0 & 0 & 0 \\ 
0 & 0 & -2\alpha ^{2}%
\end{array}%
\right) \text{ with }c=2\alpha ^{2}.
\end{equation*}

\item $\alpha =\gamma >0$ and $\beta =0$. In this case we have $\mathrm{D}=0$
and $c=0.$

\item $\alpha =\beta =\gamma =0$. In this case we have 
\begin{equation*}
\mathrm{D=}\left( 
\begin{array}{ccc}
\lambda _{1}^{1} & 0 & 0 \\ 
0 & \lambda _{1}^{1} & 0 \\ 
0 & 0 & \lambda _{1}^{1}%
\end{array}%
\right) \text{ with }c=-\lambda _{1}^{1}.
\end{equation*}
\end{itemize}
\end{theorem}

\textbf{Proof. }First assume that $\alpha \neq 0$ and $\beta =\gamma =0,$ in
this case the algebraic Ricci Soliton condition 
\eqref{Soliton}
on $G_{3}$ reduces to the following system:%
\begin{eqnarray}
&&%
\begin{array}{c}
c+\lambda _{2}^{2}+\lambda _{3}^{3}=-\left( c+\lambda _{2}^{2}\right)
=-\left( c+\lambda _{3}^{3}\right) =-\frac{\alpha ^{2}}{2},\text{ \ }%
\end{array}
\label{Sol3a} \\
&&%
\begin{array}{c}
\lambda _{2}^{1}=\lambda _{2}^{3}=\lambda _{3}^{1}=\lambda _{3}^{2}=0.%
\end{array}
\notag
\end{eqnarray}

Thus, the first equation in 
\eqref{Sol3a}
gives%
\begin{equation*}
\lambda _{2}^{2}=\lambda _{3}^{3}=-\alpha ^{2}\text{ and }c=\frac{3\alpha
^{2}}{2}.
\end{equation*}

Next, assume that $\gamma \neq 0$ and $\beta =\gamma =0$ thus 
\eqref{Soliton}
becomes%
\begin{eqnarray}
&&%
\begin{array}{c}
c+\lambda _{3}^{3}-\lambda _{2}^{2}=c+\lambda _{2}^{2}=-\left( c+\lambda
_{3}^{3}\right) =\frac{\gamma ^{2}}{2},\text{ \ }%
\end{array}
\label{Sol3b} \\
&&%
\begin{array}{c}
\lambda _{1}^{2}=\lambda _{1}^{3}=\lambda _{2}^{1}=\lambda _{2}^{3}=0.%
\end{array}
\notag
\end{eqnarray}

Therefore, from the first equation in 
\eqref{Sol3b}%
, we get%
\begin{equation*}
\lambda _{2}^{2}=-\gamma ^{2},\text{ \ }\lambda _{3}^{3}=-\gamma ^{2}\text{
\ and }c=\frac{3\gamma ^{2}}{2}.
\end{equation*}

Now, suppose that $\alpha \beta \gamma \neq 0,$ hence 
\eqref{Soliton}
reduces to%
\begin{eqnarray}
&&%
\begin{array}{c}
c=-a_{1}a_{2}-a_{1}a_{3}-\beta a_{2}-\gamma a_{3},%
\end{array}
\label{Sol3c} \\
&&%
\begin{array}{c}
c=-a_{1}a_{3}+a_{2}a_{3}+\alpha a_{1}-\beta a_{3},%
\end{array}
\notag \\
&&%
\begin{array}{c}
c=a_{2}a_{3}-a_{1}a_{2}+\alpha a_{1}-\gamma a_{3},%
\end{array}
\notag \\
&&%
\begin{array}{c}
\lambda _{1}^{2}=\lambda _{1}^{3}=\lambda _{2}^{3}=0.%
\end{array}
\notag
\end{eqnarray}

Thus, it is easy to prove that 
\eqref{Sol3c}
does not occur for all possible values of $\alpha ,\beta $ and $\gamma .$ In
fact, from the three first equations in 
\eqref{Sol3c}
it follows that 
\eqref{Sol3c}
admits solutions if and only if $\alpha =\beta =\gamma $ with $c=-\frac{%
\alpha ^{2}}{2}.$

If $\alpha \beta \neq 0$ and $\gamma =0.$ In this case, 
\eqref{Soliton}
becomes%
\begin{eqnarray}
&&%
\begin{array}{c}
c+\lambda _{1}^{1}=-\frac{1}{2}\left( \alpha ^{2}-\beta ^{2}\right) =\frac{1%
}{2}\left( \alpha ^{2}-\beta ^{2}\right) ,%
\end{array}
\label{Sol3d} \\
&&%
\begin{array}{c}
c=\frac{1}{2}\left( \alpha -\beta \right) ^{2},%
\end{array}
\notag \\
&&%
\begin{array}{c}
\lambda _{1}^{2}=\lambda _{3}^{1}=\lambda _{3}^{2}=0.%
\end{array}
\notag
\end{eqnarray}

Hence 
\eqref{Sol3d}
admits solutions if and only if $\alpha =\pm \beta .$ Note that, when $%
\alpha >0$ and $\beta >0$ we have $\alpha =\beta $ and $\lambda
_{1}^{1}=c=0. $ However, if $\alpha >0$ and $\beta <0$ we get $\lambda
_{1}^{1}=-c=-2$ $\alpha ^{2}.$

Finally, if $\alpha =\beta =\gamma =0,$ in this case the algebraic Ricci
Soliton condition 
\eqref{Soliton}
on $G_{3}$ is satisfied if and only if 
\begin{eqnarray*}
&&%
\begin{array}{c}
\begin{array}{c}
\lambda _{1}^{1}=\lambda _{2}^{2}=\lambda _{3}^{3}=-c,%
\end{array}%
\end{array}
\\
&&%
\begin{array}{c}
\lambda _{1}^{2}=\lambda _{1}^{3}=\lambda _{2}^{1}=\lambda _{2}^{3}=\lambda
_{3}^{1}=\lambda _{3}^{2}=0. \hfill \square%
\end{array}%
\end{eqnarray*}

\begin{remark}
Let $D$ be a derivation of $\mathfrak{g}_{3}$ for which $G_{3}$ is an algebraic Ricci Soliton Lie group. 
Hence, standard calculations prove that $D\neq \mathrm{ad}X$, for any $X\in \mathfrak{g}_{3}$.
\end{remark}

\begin{remark}
Above soliton have been already constructed.

\begin{itemize}
\item $\alpha >0$ and $\beta =\gamma =0.$ This is $(H_{3},g_{1})$ (see \cite%
{O10} and \cite{O11}).

\item $\gamma <0$ and $\alpha =\beta =0.$ This is $(H_{3},g_{2})$ (see \cite%
{O11}).

\item $\alpha =\beta =\gamma \neq 0$. This is a negative constant metric
(see \cite{N79}).

\item $\alpha =\beta >0$ and $\gamma =0$. This is a flat metric on $E(2)$
(see \cite{N79}).

\item $\beta =-\alpha <0$ and $\gamma =0$. In particular, $\alpha =1$. This
is $(E(1,1),g_{1})$ (see \cite{O10} and \cite{O11}).

\item $\gamma =-\alpha <0$ and $\beta =0$. In particular, $\alpha =1$. This
is $(E(2),g_{1})$ (see \cite{O10} and \cite{O11}).

\item $\alpha =\gamma >0$ and $\beta =0$. This is a flat metric on $E(1,1)$
(see \cite{N79}).

\item $\alpha =\beta =\gamma =0$. This is a Gaussian soliton on $%
\mathbb{R}
^{3}$ (see \cite{P02}).
\end{itemize}
\end{remark}

\subsection{Algebraic Ricci Solitons of $G_{4}$}
Finally consider the Lorentzian Lie group $G_{4},$ the following lemma was
proven in \cite{C07a}.

\begin{lemma}
Let $\{e_{1},e_{2},e_{3}\}$ be the pseudo-orthonormal basis used in 
\eqref{g4}%
. Then%
\begin{eqnarray*}
\nabla _{e_{1}}e_{1} &=&0,\text{ \ \ \ \ \ \ \ \ }\nabla
_{e_{2}}e_{1}=e_{2}+b_{2}e_{3},\text{ \ \ \ \ \ }\nabla
_{e_{3}}e_{1}=b_{3}e_{2}-e_{3}, \\
\nabla _{e_{1}}e_{2} &=&b_{1}e_{3},\text{ \ \ \ \ }\nabla
_{e_{2}}e_{2}=-e_{1},\text{\ \ \ \ \ \ \ \ \ \ \ \ }\nabla
_{e_{3}}e_{2}=-b_{3}e_{1}, \\
\nabla _{e_{1}}e_{3} &=&b_{1}e_{2},\text{ \ \ \ \ }\nabla
_{e_{2}}e_{3}=b_{2}e_{1},\text{ \ \ \ \ \ \ \ \ \ \ \ }\nabla
_{e_{3}}e_{3}=-e_{1},
\end{eqnarray*}%
where we put%
\begin{equation*}
b_{1}=\frac{\alpha }{2}+\eta -\beta ,\text{ \ \ \ }b_{2}=\frac{\alpha }{2}%
-\eta ,\text{ \ \ \ }b_{3}=\frac{\alpha }{2}+\eta .
\end{equation*}%
The not always null components of the curvature tensor are given by%
\begin{eqnarray*}
R_{1212} &=&\left( 2\eta -\beta \right) b_{3}-b_{1}b_{2}-1,\text{ \ }%
R_{1313}=b_{1}b_{3}+\beta b_{2}+1,\text{ \ \ } \\
R_{2323} &=&-\left( b_{2}b_{3}+\alpha b_{1}+1\right) ,\text{ \ \ \ \ }%
R_{1213}=2\eta -\beta +b_{1}+b_{2}.
\end{eqnarray*}%
The Ricci operator is given by%
\begin{equation*}
\mathrm{Ric}=\left( 
\begin{array}{ccc}
-\frac{\alpha ^{2}}{2} & 0 & 0 \\ 
0 & \frac{\alpha ^{2}}{2}+2\eta \left( \alpha -\beta \right) -\alpha \beta +2
& \alpha +2\eta -2\beta \\ 
0 & -\alpha -2\eta +2\beta & \frac{\alpha ^{2}}{2}-\alpha \beta -2+2\eta
\beta%
\end{array}%
\right) .
\end{equation*}
\end{lemma}

Now, put $\mathrm{D}e_{l}=\lambda _{l}^{1}e_{1}+\lambda
_{l}^{2}e_{2}+\lambda _{l}^{3}e_{3}$ where $\{e_{1},e_{2},e_{3}\}$ is the
pseudo-orthonormal basis used in 
\eqref{g4}%
. Hence, $\mathrm{D}\in \mathrm{Der}\left( \mathfrak{g}_{4}\right) $ if and
only if%
\begin{eqnarray}
&&%
\begin{array}{c}
\alpha \lambda _{1}^{3}+\left( 2\eta -\beta \right) \lambda _{3}^{1}-\lambda
_{2}^{1}=0,%
\end{array}
\label{sysD4} \\
&&%
\begin{array}{c}
\beta \lambda _{2}^{3}+\left( 2\eta -\beta \right) \lambda _{3}^{2}+\lambda
_{1}^{1}=0,%
\end{array}
\notag \\
&&%
\begin{array}{c}
\left( 2\eta -\beta \right) \left( \lambda _{1}^{1}+\lambda _{2}^{2}-\lambda
_{3}^{3}\right) +2\lambda _{2}^{3}=0,%
\end{array}
\notag \\
&&%
\begin{array}{c}
\alpha \lambda _{1}^{2}+\beta \lambda _{2}^{1}-\lambda _{3}^{1}=0,%
\end{array}
\notag \\
&&%
\begin{array}{c}
\beta \left( \lambda _{1}^{1}-\lambda _{2}^{2}+\lambda _{3}^{3}\right)
+2\lambda _{3}^{2}=0,%
\end{array}
\notag \\
&&%
\begin{array}{c}
\alpha \left( -\lambda _{1}^{1}+\lambda _{2}^{2}+\lambda _{3}^{3}\right) =0.%
\end{array}
\notag
\end{eqnarray}

Again, a straightforward computation lead to the following:

\begin{lemma}
The solutions of 
\eqref{sysD4}
are the following:

\begin{itemize}
\item If $\alpha =0$ and $\beta =\eta :$%
\begin{equation*}
\mathrm{D=}\left( 
\begin{array}{ccc}
\lambda _{1}^{1} & \lambda _{2}^{1} & \eta \lambda _{2}^{1} \\ 
\lambda _{1}^{2} & \lambda _{2}^{2} & \frac{\eta }{2}\left( -\lambda
_{1}^{1}+\lambda _{2}^{2}-\lambda _{3}^{3}\right) \\ 
\lambda _{1}^{3} & \frac{\eta }{2}\left( -\lambda _{1}^{1}-\lambda
_{2}^{2}+\lambda _{3}^{3}\right) & \lambda _{3}^{3}%
\end{array}%
\right) .
\end{equation*}

\item If $\alpha =0$ and $\beta \neq \eta :$%
\begin{equation*}
\mathrm{D=}\left( 
\begin{array}{ccc}
0 & 0 & 0 \\ 
\lambda _{1}^{2} & \lambda _{2}^{2} & \frac{\beta }{2}\left( \lambda
_{2}^{2}-\lambda _{3}^{3}\right) \\ 
\lambda _{1}^{3} & \frac{\left( 2\eta -\beta \right) }{2}\left( -\lambda
_{2}^{2}+\lambda _{3}^{3}\right) & \lambda _{3}^{3}%
\end{array}%
\right) .
\end{equation*}

\item If $\alpha \neq 0$ and $\beta =\eta :$ 
\begin{equation*}
\mathrm{D=}\left( 
\begin{array}{ccc}
\lambda _{2}^{2}+\lambda _{3}^{3} & \eta \left( \lambda _{3}^{1}-\alpha
\lambda _{1}^{2}\right) & \lambda _{3}^{1} \\ 
\lambda _{1}^{2} & \lambda _{2}^{2} & -\eta \lambda _{3}^{3} \\ 
-\eta \lambda _{1}^{2} & -\eta \lambda _{2}^{2} & \lambda _{3}^{3}%
\end{array}%
\right) .
\end{equation*}

\item If $\alpha \neq 0$ and $\beta \neq \eta :$%
\begin{equation*}
\mathrm{D=}\left( 
\begin{array}{ccc}
0 & \lambda _{2}^{1} & \lambda _{3}^{1} \\ 
-\frac{\eta }{\alpha }\lambda _{2}^{1}+\frac{1}{\alpha }\lambda _{3}^{1} & 
\lambda _{2}^{2} & \beta \lambda _{2}^{2} \\ 
\frac{1}{\alpha }\lambda _{2}^{1}-\frac{\left( 2\eta -\beta \right) }{\alpha 
}\lambda _{3}^{1} & \left( -2\eta +\beta \right) \lambda _{2}^{2} & -\lambda
_{2}^{2}%
\end{array}%
\right) .
\end{equation*}
\end{itemize}
\end{lemma}

We now prove the following.

\begin{theorem}
Consider the three-dimensional Lorentzian Lie group $G_{4}$. Then, $G_{4}$
is an algebraic Ricci Soliton Lie group if and only if $\alpha =0$ and $%
\beta =\eta .$ In particular%
\begin{equation*}
\mathrm{D}=0\text{ and }c=0.
\end{equation*}
\end{theorem}

\textbf{Proof. }Assume that $\left( \alpha ,\beta \right) =\left( 0,\eta
\right) ,$ thus the left-invariant Lorentzian metric in this case is flat.
Therefore, one can show that 
\eqref{Soliton}
on $G_{4}$ is satisfied only for $\mathrm{D}=0$ and $c=0.$

Now suppose that $\alpha =0$ and $\beta \neq \eta ,$ then the algebraic
Ricci Soliton condition 
\eqref{Soliton}
on $G_{4}$ reduces to the following system: 
\begin{eqnarray}
&&%
\begin{array}{c}
c+\lambda _{2}^{2}=-\left( c+\lambda _{3}^{3}\right) =2\left( 1-\eta \beta
\right) ,%
\end{array}
\label{Sol4a} \\
&&%
\begin{array}{c}
\frac{\beta }{2}\left( \lambda _{2}^{2}-\lambda _{3}^{3}\right) =\frac{%
\left( 2\eta -\beta \right) }{2}\left( \lambda _{2}^{2}-\lambda
_{3}^{3}\right) =2\left( \eta -\beta \right) ,%
\end{array}
\notag \\
&&%
\begin{array}{c}
\lambda _{1}^{2}=\lambda _{1}^{3}=c=0.%
\end{array}
\notag
\end{eqnarray}

Thus, the first equation in 
\eqref{Sol4a}
gives%
\begin{equation*}
\lambda _{2}^{2}=-\lambda _{3}^{3}=2\left( 1-\eta \beta \right) .
\end{equation*}

Replacing $\lambda _{2}^{2}$ and $\lambda _{3}^{3}$ in the second equation
of 
\eqref{Sol4a}%
, we get%
\begin{equation*}
\left( \eta -\beta \right) \left( 1-\eta \beta \right) =0,
\end{equation*}%
which is a contradiction, since $\beta \neq \eta =\pm 1.$

If $\alpha \neq 0$ and $\beta =\eta $ we find that 
\eqref{Soliton}
is equivalent to the following system:%
\begin{eqnarray}
&&%
\begin{array}{c}
c+\lambda _{2}^{2}+\lambda _{3}^{3}=-\frac{\alpha ^{2}}{2},%
\end{array}
\label{Sol4b} \\
&&%
\begin{array}{c}
c+\lambda _{2}^{2}=\frac{\alpha ^{2}}{2}+\eta \alpha ,%
\end{array}
\notag \\
&&%
\begin{array}{c}
c+\lambda _{3}^{3}=\frac{\alpha ^{2}}{2}-\eta \alpha ,%
\end{array}
\notag \\
&&%
\begin{array}{c}
\lambda _{2}^{2}=-\lambda _{3}^{3}=\eta \alpha ,%
\end{array}
\notag \\
&&%
\begin{array}{c}
\lambda _{1}^{2}=\lambda _{3}^{1}=0.%
\end{array}
\notag
\end{eqnarray}

Hence, it follows, after replacing $\lambda _{2}^{2}$ and $\lambda _{3}^{3}$
in the first equation of 
\eqref{Sol4b}%
, that $\alpha =0$ which is a contradiction.

For $\alpha \neq 0$ and $\beta \neq \eta ,$ system 
\eqref{Soliton}
is satisfied if and only if%
\begin{eqnarray}
&&%
\begin{array}{c}
c=-\frac{\alpha ^{2}}{2},%
\end{array}
\label{Sol4c} \\
&&%
\begin{array}{c}
c+\lambda _{2}^{2}=\frac{\alpha ^{2}}{2}+2\eta \left( \alpha -\beta \right)
-\alpha \beta +2,%
\end{array}
\notag \\
&&%
\begin{array}{c}
c-\lambda _{2}^{2}=\frac{\alpha ^{2}}{2}-\alpha \beta -2+2\eta \beta ,%
\end{array}
\notag \\
&&%
\begin{array}{c}
\left( -2\eta +\beta \right) \lambda _{2}^{2}=-\alpha -2\eta +2\beta%
\end{array}
\notag \\
&&%
\begin{array}{c}
\beta \lambda _{2}^{2}=\alpha +2\eta -2\beta%
\end{array}
\notag \\
&&%
\begin{array}{c}
\lambda _{2}^{1}=\lambda _{3}^{1}=0.%
\end{array}
\notag
\end{eqnarray}

The fourth and fifth equation in 
\eqref{Sol4c}
implies $\lambda _{2}^{2}=0$ and $\alpha =2\left( \beta -\eta \right) $.
Then, by the second equation of 
\eqref{Sol4c}%
, we find that $\left( \beta -\eta \right) =0$ and so, this case does not
occur.

\begin{remark}
If $\alpha =0$ and $\beta =\eta $ , this case correspond to the flat metric
on the Heisenberg group $H_{3}$, obtained by Nomizu (see \cite{N79}).
\end{remark}

\begin{remark}
\label{maintheorem2b} Put $Y=\sum Y^{i}e_{i}$ where $\{e_{1},e_{2},e_{3}\}$
denotes the basis used in 
\eqref{g4}
and let $g$ be the left-invariant Lorentzian metric for which $%
\{e_{1},e_{2},e_{3}\}$ is the pseudo-orthonormal for $G_{4}$. Therefore%
\begin{equation*}
L_{Y}g=\left( 
\begin{array}{ccc}
0 & (\alpha -\beta )Y^{3}-Y^{2} & -(\alpha +2\eta -\beta )Y^{2}-Y^{3} \\ 
(\alpha -\beta )Y^{3}-Y^{2} & 2Y^{1} & 2\eta Y^{1} \\ 
-(\alpha +2\eta -\beta )Y^{2}-Y^{3} & 2\eta Y^{1} & 2Y^{1}%
\end{array}%
\right) ,
\end{equation*}%
and hence, we have a Ricci soliton equation 
\eqref{Ricci soliton}
is satisfied if and only if%
\begin{eqnarray*}
&&%
\begin{array}{c}
\dfrac{\alpha ^{2}}{2}+2\eta (\alpha -\beta )-\alpha \beta +2=2Y^{1}+c,%
\end{array}
\\
&&%
\begin{array}{c}
-\dfrac{\alpha ^{2}}{2}+\alpha \beta +2-2\eta \beta =2Y^{1}-c,%
\end{array}
\\
&&%
\begin{array}{c}
(\alpha +2\eta -\beta )Y^{2}+Y^{3}=0,%
\end{array}
\\
&&%
\begin{array}{c}
(\alpha -\beta )Y^{3}-Y^{2}=0,%
\end{array}
\\
&&%
\begin{array}{c}
\alpha +2\eta -2\beta =2\eta Y^{1},%
\end{array}
\\
&&%
\begin{array}{c}
c=-\dfrac{\alpha ^{2}}{2}.%
\end{array}%
\end{eqnarray*}%
Standard calculation show that, whenever $\left( \alpha ,\beta \right)
=\left( 0,\eta \right) $, we have $Y^{1}=0$ and so, the Ricci tensor in this
case vanish identically, in this case the metric is flat. On the other hand,
if $\alpha =0$ and $\beta \neq \eta $ then $Y^{2}=Y^{3}=0$ and $c=0,$ thus $%
Y=(-\eta \beta +1)e_{1}$ this is a nonflat steady Ricci soliton on $E(1,1)$.
Now, assume that $\alpha \neq 0$, then $Y=-\dfrac{\eta }{2}\alpha
e_{1}+Y^{2}\left( e_{2}-\eta e_{3}\right) $ and $c=-\dfrac{\alpha ^{2}}{2}.$
Therefore we obtain an expanding Ricci soliton on $SL(2,\mathbb{R})$ which
is not an algebraic soliton, since $\alpha \neq 0$ (see the previous Theorem).
\end{remark}

\section{\protect\LARGE Algebraic Ricci Solitons of 3-dimensional
non-unimodular Lie groups}

\setcounter{equation}{0}

Now, we shall inspect our three basic, non-unimodular, algebras case by case.

\subsection{Algebraic Ricci Solitons of $G_{5}$}
We start by recalling the following result on the curvature tensor and the
Ricci operator of a three-dimensional Lorentzian Lie group $G_{5}.$

\begin{lemma}[\protect\cite{C07a}]
Let $\{e_{1},e_{2},e_{3}\}$ be the pseudo-orthonormal basis used in 
\eqref{g5}%
. Then%
\begin{eqnarray*}
\nabla _{e_{1}}e_{1} &=&\alpha e_{3},\text{ \ \ \ \ \ \ \ \ \ \ \ \ \ \ \ \ }%
\nabla _{e_{2}}e_{1}=\frac{\beta +\gamma }{2}e_{3},\text{ \ \ \ \ \ \ \ \ \
\ \ \ }\nabla _{e_{3}}e_{1}=-\frac{\beta -\gamma }{2}e_{2}, \\
\nabla _{e_{1}}e_{2} &=&\frac{\beta +\gamma }{2}e_{3},\text{ \ \ \ \ \ \ \ \
\ \ \ }\nabla _{e_{2}}e_{2}=\delta e_{3},\text{ \ \ \ \ \ \ \ \ \ \ \ \ \ \
\ \ \ \ }\nabla _{e_{3}}e_{2}=\frac{\beta -\gamma }{2}e_{1}, \\
\nabla _{e_{1}}e_{3} &=&\alpha e_{1}+\frac{\beta +\gamma }{2}e_{2},\text{ \
\ \ }\nabla _{e_{2}}e_{3}=\frac{\beta +\gamma }{2}e_{1}+\delta e_{2},\text{
\ \ \ \ \ }\nabla _{e_{3}}e_{3}=0,
\end{eqnarray*}%
and the possible non-zero components of the curvature tensor are given by%
\begin{eqnarray*}
R_{1212} &=&\alpha \delta -\frac{\left( \beta +\gamma \right) ^{2}}{4},\text{
\ \ \ }R_{1313}=-\alpha ^{2}-\frac{\beta \left( \beta +\gamma \right) }{2}-%
\frac{\beta ^{2}-\gamma ^{2}}{4},\text{ \ \ \ } \\
R_{2323} &=&-\delta ^{2}-\frac{\gamma \left( \beta +\gamma \right) }{2}+%
\frac{\beta ^{2}-\gamma ^{2}}{4}.
\end{eqnarray*}%
The Ricci operator is given by%
\begin{equation*}
\mathrm{Ric}=\left( 
\begin{array}{ccc}
\alpha ^{2}+\alpha \delta +\frac{\beta ^{2}-\gamma ^{2}}{2} & 0 & 0 \\ 
0 & \alpha \delta +\delta ^{2}-\frac{\beta ^{2}-\gamma ^{2}}{2} & 0 \\ 
0 & 0 & \alpha ^{2}+\delta ^{2}+\frac{\left( \beta +\gamma \right) ^{2}}{2}%
\end{array}%
\right) .
\end{equation*}
\end{lemma}

Next, put $\mathrm{D}e_{l}=\lambda _{l}^{1}e_{1}+\lambda
_{l}^{2}e_{2}+\lambda _{l}^{3}e_{3}$ where $\{e_{1},e_{2},e_{3}\}$ is the
pseudo-orthonormal basis used in 
\eqref{g5}%
. Then, $\mathrm{D}\in \mathrm{Der}\left( \mathfrak{g}_{5}\right) $ if and
only if%
\begin{eqnarray}
&&%
\begin{array}{c}
\alpha \lambda _{2}^{3}-\gamma \lambda _{1}^{3}=0,%
\end{array}
\label{sysD5} \\
&&%
\begin{array}{c}
\beta \lambda _{2}^{3}-\delta \lambda _{1}^{3}=0,%
\end{array}
\notag \\
&&%
\begin{array}{c}
\alpha \lambda _{3}^{3}+\gamma \lambda _{1}^{2}-\beta \lambda _{2}^{1}=0,%
\end{array}
\notag \\
&&%
\begin{array}{c}
\beta \left( \lambda _{1}^{1}-\lambda _{2}^{2}+\lambda _{3}^{3}\right)
+\left( \delta -\alpha \right) \lambda _{1}^{2}=0,%
\end{array}
\notag \\
&&%
\begin{array}{c}
\alpha \lambda _{1}^{3}+\beta \lambda _{2}^{3}=0,%
\end{array}
\notag \\
&&%
\begin{array}{c}
\gamma \left( -\lambda _{1}^{1}+\lambda _{2}^{2}+\lambda _{3}^{3}\right)
+\left( \alpha -\delta \right) \lambda _{2}^{1}=0,%
\end{array}
\notag \\
&&%
\begin{array}{c}
\delta \lambda _{3}^{3}+\beta \lambda _{2}^{1}-\gamma \lambda _{1}^{2}=0,%
\end{array}
\notag \\
&&%
\begin{array}{c}
\gamma \lambda _{1}^{3}+\beta \lambda _{2}^{3}=0.%
\end{array}
\notag
\end{eqnarray}%
In determining the solutions of 
\eqref{sysD5}
we must also take into account that, by 
\eqref{g5}%
, $\alpha +\delta \neq 0$ and $\alpha \gamma +\beta \delta =0.$ A
straightforward computation lead to the following:

\begin{lemma}
The solutions of 
\eqref{sysD5}
are the following:

\begin{itemize}
\item If $\left( \beta ,\gamma \right) \neq \left( 0,0\right) :$%
\begin{equation*}
\mathrm{D=}\left( 
\begin{array}{ccc}
\lambda _{1}^{1} & \gamma k & \lambda _{3}^{1} \\ 
\beta k & \lambda _{1}^{1}+\left( \delta -\alpha \right) k & \lambda _{3}^{2}
\\ 
0 & 0 & 0%
\end{array}%
\right) \text{ with }k\in 
\mathbb{R}
.
\end{equation*}

\item If $\left( \beta ,\gamma \right) =\left( 0,0\right) $ and $\delta \neq
\alpha :$%
\begin{equation*}
\mathrm{D=}\left( 
\begin{array}{ccc}
\lambda _{1}^{1} & 0 & \lambda _{3}^{1} \\ 
0 & \lambda _{2}^{2} & \lambda _{3}^{2} \\ 
0 & 0 & 0%
\end{array}%
\right) .
\end{equation*}

\item If $\left( \beta ,\gamma \right) =\left( 0,0\right) $ and $\delta
=\alpha :$%
\begin{equation*}
\mathrm{D=}\left( 
\begin{array}{ccc}
\lambda _{1}^{1} & \lambda _{2}^{1} & \lambda _{3}^{1} \\ 
\lambda _{1}^{2} & \lambda _{2}^{2} & \lambda _{3}^{2} \\ 
0 & 0 & 0%
\end{array}%
\right) .
\end{equation*}
\end{itemize}
\end{lemma}

Using the above lemma we prove the following:

\begin{theorem}
Consider the three-dimensional Lorentzian Lie group $G_{5}$. Then, $G_{5}$
is an algebraic Ricci Soliton Lie group if and only if one of the following
statements holds true:

\begin{itemize}
\item $\left( \beta ,\gamma \right) \neq \left( 0,0\right) $ and $\alpha
^{2}+\beta ^{2}=\gamma ^{2}+\delta ^{2}.$ In this case we have%
\begin{equation*}
\mathrm{D=}\left( 
\begin{array}{ccc}
\alpha \delta -\beta \gamma -\alpha ^{2}-\beta ^{2} & 0 & 0 \\ 
0 & \alpha \delta -\beta \gamma -\alpha ^{2}-\beta ^{2} & 0 \\ 
0 & 0 & 0%
\end{array}%
\right) \text{and }c=\alpha ^{2}+\delta ^{2}+\frac{\left( \beta +\gamma
\right) ^{2}}{2}.
\end{equation*}

\item $\left( \beta ,\gamma \right) =\left( 0,0\right) .$ In particular%
\begin{equation*}
\mathrm{D=}\left( 
\begin{array}{ccc}
\delta \left( \alpha -\delta \right) & 0 & 0 \\ 
0 & \alpha \left( \delta -\alpha \right) & 0 \\ 
0 & 0 & 0%
\end{array}%
\right) \text{ and }c=\alpha ^{2}+\delta ^{2}.
\end{equation*}
\end{itemize}
\end{theorem}

\textbf{Proof. }Assume that $\left( \beta ,\gamma \right) \neq \left(
0,0\right) ,$ we find that 
\eqref{Soliton}
is equivalent to the following system:%
\begin{eqnarray}
&&%
\begin{array}{c}
c+\lambda _{1}^{1}=\alpha ^{2}+\alpha \delta +\frac{\beta ^{2}-\gamma ^{2}}{2%
},%
\end{array}
\label{Sol5a} \\
&&%
\begin{array}{c}
c+\lambda _{1}^{1}+\left( \delta -\alpha \right) =\alpha \delta +\delta ^{2}-%
\frac{\beta ^{2}-\gamma ^{2}}{2},%
\end{array}
\notag \\
&&%
\begin{array}{c}
c=\alpha ^{2}+\delta ^{2}+\frac{\left( \beta +\gamma \right) ^{2}}{2},%
\end{array}
\notag \\
&&%
\begin{array}{c}
\lambda _{3}^{1}=\lambda _{3}^{2}=0,%
\end{array}
\notag \\
&&%
\begin{array}{c}
\beta k=\gamma k=0.%
\end{array}
\notag
\end{eqnarray}

Since $\left( \beta ,\gamma \right) \neq \left( 0,0\right) $, we have, from
the last equation of 
\eqref{Sol5a}%
, that $k=0.$ Therefore, we get 
\begin{equation*}
\lambda _{1}^{1}=\alpha \delta -\beta \gamma -\alpha ^{2}-\beta ^{2}\text{
and }c=\alpha ^{2}+\delta ^{2}+\frac{\left( \beta +\gamma \right) ^{2}}{2},
\end{equation*}%
but only if $\alpha ^{2}+\beta ^{2}=\gamma ^{2}+\delta ^{2}.$ Assume now
that, $\left( \beta ,\gamma \right) =\left( 0,0\right) $ and $\delta \neq
\alpha $. In this case, the algebraic Ricci Soliton condition 
\eqref{Soliton}
on $G_{5}$ reduces to the following system:%
\begin{eqnarray*}
&&%
\begin{array}{c}
c+\lambda _{1}^{1}=\alpha ^{2}+\alpha \delta ,%
\end{array}
\\
&&%
\begin{array}{c}
c+\lambda _{2}^{2}=\alpha \delta +\delta ^{2},%
\end{array}
\\
&&%
\begin{array}{c}
c=\alpha ^{2}+\delta ^{2},%
\end{array}
\\
&&%
\begin{array}{c}
\lambda _{3}^{1}=\lambda _{3}^{2}=0.%
\end{array}%
\end{eqnarray*}%
Thus, 
\begin{equation*}
\lambda _{1}^{1}=\delta \left( \alpha -\delta \right) \text{ and }\lambda
_{2}^{2}=\alpha \left( \delta -\alpha \right) .
\end{equation*}

In the remaining case, $\left( \beta ,\gamma \right) =\left( 0,0\right) $
and $\delta =\alpha $ and so, the algebraic Ricci Soliton condition 
\eqref{Soliton}
on $G_{5}$ is satisfied if and only if%
\begin{equation}
\lambda _{1}^{1}=\lambda _{2}^{2}=\lambda _{1}^{2}=\lambda _{2}^{1}=\lambda
_{3}^{1}=\lambda _{3}^{2}=0\text{ and }c=2\alpha ^{2}.  \notag
\end{equation}

\begin{remark} 
Let $D$ denote the derivation of $\mathfrak{g}_{5}$ given in the above Theorem. 
Then, it is easy to prove that $D$ is an outer derivation. 
\end{remark}


\subsection{Algebraic Ricci Solitons of $G_{6}$}
Next consider the three-dimensional Lorentzian Lie group $G_{6}$. We recall
the following result :

\begin{lemma}[\protect\cite{C07a}]
Let $\{e_{1},e_{2},e_{3}\}$ be the pseudo-orthonormal basis used in 
\eqref{g6}%
. Then%
\begin{eqnarray*}
\nabla _{e_{1}}e_{1} &=&0,\text{ \ \ \ \ \ \ \ \ \ \ \ \ \ \ \ }\nabla
_{e_{2}}e_{1}=-\alpha e_{2}-\frac{\beta -\gamma }{2}e_{3},\text{ \ \ \ \ \ \
\ \ }\nabla _{e_{3}}e_{1}=\frac{\beta -\gamma }{2}e_{2}-\delta e_{3}, \\
\nabla _{e_{1}}e_{2} &=&\frac{\beta +\gamma }{2}e_{3},\text{ \ \ \ \ \ \ }%
\nabla _{e_{2}}e_{2}=\alpha e_{1},\text{ \ \ \ \ \ \ \ \ \ \ \ \ \ \ \ \ \ \
\ \ \ \ \ \ }\nabla _{e_{3}}e_{2}=-\frac{\beta -\gamma }{2}e_{1}, \\
\nabla _{e_{1}}e_{3} &=&\frac{\beta +\gamma }{2}e_{2},\text{ \ \ \ \ \ \ }%
\nabla _{e_{2}}e_{3}=-\frac{\beta -\gamma }{2}e_{1},\text{ \ \ \ \ \ \ \ \ \
\ \ \ \ \ \ \ }\nabla _{e_{3}}e_{3}=-\delta e_{1},
\end{eqnarray*}%
and the only possibly non-vanishing components of the curvature tensor are
given by%
\begin{eqnarray*}
R_{1212} &=&-\alpha ^{2}+\frac{\beta ^{2}-\gamma ^{2}}{4}+\frac{\beta \left(
\beta -\gamma \right) }{2},\text{ \ \ }R_{1313}=\delta ^{2}+\frac{\beta
^{2}-\gamma ^{2}}{4}+\frac{\gamma \left( \beta -\gamma \right) }{2},\text{ }
\\
R_{2323} &=&\alpha \delta +\frac{\left( \beta -\gamma \right) ^{2}}{4}.
\end{eqnarray*}%
The Ricci operator is given by%
\begin{equation*}
\mathrm{Ric}=\left( 
\begin{array}{ccc}
-\alpha ^{2}-\delta ^{2}+\frac{\left( \beta -\gamma \right) ^{2}}{2} & 0 & 0
\\ 
0 & -\alpha ^{2}-\alpha \delta +\frac{\beta ^{2}-\gamma ^{2}}{2} & 0 \\ 
0 & 0 & -\delta ^{2}-\alpha \delta -\frac{\beta ^{2}-\gamma ^{2}}{2}%
\end{array}%
\right) .
\end{equation*}
\end{lemma}

Now, put $\mathrm{D}e_{l}=\lambda _{l}^{1}e_{1}+\lambda
_{l}^{2}e_{2}+\lambda _{l}^{3}e_{3}$ an endomorphism of $\mathfrak{g}_{6}$
where $\{e_{1},e_{2},e_{3}\}$ is the pseudo-orthonormal basis used in 
\eqref{g6}%
. Then, $\mathrm{D}\in \mathrm{Der}\left( \mathfrak{g}_{6}\right) $ if and
only if%
\begin{eqnarray}
&&%
\begin{array}{c}
\alpha \lambda _{2}^{1}+\beta \lambda _{3}^{1}=0,%
\end{array}
\label{sysD6} \\
&&%
\begin{array}{c}
\alpha \lambda _{1}^{1}+\gamma \lambda _{2}^{3}-\beta \lambda _{3}^{2}=0,%
\end{array}
\notag \\
&&%
\begin{array}{c}
\beta \left( \lambda _{1}^{1}+\lambda _{2}^{2}-\lambda _{3}^{3}\right)
+\left( \delta -\alpha \right) \lambda _{2}^{3}=0,%
\end{array}
\notag \\
&&%
\begin{array}{c}
\gamma \lambda _{2}^{1}+\delta \lambda _{3}^{1}=0,%
\end{array}
\notag \\
&&%
\begin{array}{c}
\gamma \left( \lambda _{1}^{1}-\lambda _{2}^{2}+\lambda _{3}^{3}\right)
+\left( \alpha -\delta \right) \lambda _{3}^{2}=0,%
\end{array}
\notag \\
&&%
\begin{array}{c}
\delta \lambda _{1}^{1}+\beta \lambda _{3}^{2}-\gamma \lambda _{2}^{3}=0,%
\end{array}
\notag \\
&&%
\begin{array}{c}
\gamma \lambda _{2}^{1}-\alpha \lambda _{3}^{1}=0,%
\end{array}
\notag \\
&&%
\begin{array}{c}
\delta \lambda _{2}^{1}-\beta \lambda _{3}^{1}=0.%
\end{array}
\notag
\end{eqnarray}%
In determining the solutions of 
\eqref{sysD6}
we must also take into account that, by 
\eqref{g6}%
, $\alpha +\delta \neq 0$ and $\alpha \gamma -\beta \delta =0.$ Routine but
long calculation lead to the following:

\begin{lemma}
The solutions of 
\eqref{sysD6}
are the following:

\begin{itemize}
\item If $\left( \beta ,\gamma \right) \neq \left( 0,0\right) :$%
\begin{equation*}
\mathrm{D=}\left( 
\begin{array}{ccc}
0 & 0 & 0 \\ 
\lambda _{1}^{2} & \lambda _{2}^{2} & \gamma k \\ 
\lambda _{1}^{3} & \beta k & \lambda _{2}^{2}+\left( \delta -\alpha \right) k%
\end{array}%
\right) \text{ with }k\in 
\mathbb{R}
.
\end{equation*}

\item If $\left( \beta ,\gamma \right) =\left( 0,0\right) $ and $\delta \neq
\alpha :$%
\begin{equation*}
\mathrm{D=}\left( 
\begin{array}{ccc}
0 & 0 & 0 \\ 
\lambda _{1}^{2} & \lambda _{2}^{2} & 0 \\ 
\lambda _{1}^{3} & 0 & \lambda _{3}^{3}%
\end{array}%
\right) .
\end{equation*}

\item If $\left( \beta ,\gamma \right) =\left( 0,0\right) $ and $\delta
=\alpha :$%
\begin{equation*}
\mathrm{D=}\left( 
\begin{array}{ccc}
0 & 0 & 0 \\ 
\lambda _{1}^{2} & \lambda _{2}^{2} & \lambda _{3}^{2} \\ 
\lambda _{1}^{3} & \lambda _{2}^{3} & \lambda _{3}^{3}%
\end{array}%
\right) .
\end{equation*}
\end{itemize}
\end{lemma}

Using this result, we now prove the following.

\begin{theorem}
Consider the three-dimensional Lorentzian Lie group $G_{6}$. Then, $G_{6}$
is an algebraic Ricci Soliton Lie group if and only if one of the following
statements holds true:

\begin{itemize}
\item $\left( \beta ,\gamma \right) \neq \left( 0,0\right) $ and $\alpha
^{2}-\beta ^{2}=\delta ^{2}-\gamma ^{2}.$ In this case we have%
\begin{equation*}
\mathrm{D=}\left( 
\begin{array}{ccc}
0 & 0 & 0 \\ 
0 & \alpha ^{2}-\beta ^{2}+\beta \gamma -\alpha \delta & 0 \\ 
0 & 0 & \alpha ^{2}-\beta ^{2}+\beta \gamma -\alpha \delta%
\end{array}%
\right) \text{and }c=-\alpha ^{2}-\delta ^{2}+\frac{\left( \beta -\gamma
\right) ^{2}}{2}.
\end{equation*}

\item $\left( \beta ,\gamma \right) =\left( 0,0\right) .$ In particular%
\begin{equation*}
\mathrm{D=}\left( 
\begin{array}{ccc}
0 & 0 & 0 \\ 
0 & \delta \left( \delta -\alpha \right) & 0 \\ 
0 & 0 & \alpha \left( \alpha -\delta \right)%
\end{array}%
\right) \text{ and }c=-\left( \alpha ^{2}+\delta ^{2}\right) .
\end{equation*}
\end{itemize}
\end{theorem}

\textbf{Proof. }If\textbf{\ }$\left( \beta ,\gamma \right) \neq \left(
0,0\right) .$ This case is quite similar to the corresponding one of $G_{5}$%
. In fact, the algebraic Ricci Soliton condition 
\eqref{Soliton}
on $G_{6}$ is satisfied if and only if%
\begin{eqnarray*}
&&%
\begin{array}{c}
c+\lambda _{2}^{2}=-\alpha ^{2}-\alpha \delta +\frac{\beta ^{2}-\gamma ^{2}}{%
2},%
\end{array}
\\
&&%
\begin{array}{c}
c+\lambda _{2}^{2}+\left( \delta -\alpha \right) k=-\delta ^{2}-\alpha
\delta -\frac{\beta ^{2}-\gamma ^{2}}{2},%
\end{array}
\\
&&%
\begin{array}{c}
c=-\alpha ^{2}-\delta ^{2}+\frac{\left( \beta -\gamma \right) ^{2}}{2},%
\end{array}
\\
&&%
\begin{array}{c}
\lambda _{1}^{2}=\lambda _{1}^{3}=0,%
\end{array}
\\
&&%
\begin{array}{c}
\beta k=\gamma k=0.%
\end{array}%
\end{eqnarray*}

We obtain, since $\left( \beta ,\gamma \right) \neq \left( 0,0\right) $,
that $k=0$ and so, 
\begin{equation*}
\lambda _{2}^{2}=\alpha ^{2}-\beta ^{2}+\beta \gamma -\alpha \delta ,
\end{equation*}%
but only if $\alpha ^{2}-\beta ^{2}=\delta ^{2}-\gamma ^{2}.$ If $\left(
\beta ,\gamma \right) \neq \left( 0,0\right) $ and $\delta \neq \alpha $,
the algebraic Ricci Soliton condition 
\eqref{Soliton}
on $G_{6}$ becomes:%
\begin{eqnarray*}
&&%
\begin{array}{c}
c+\lambda _{2}^{2}=-\alpha ^{2}-\alpha \delta ,%
\end{array}
\\
&&%
\begin{array}{c}
c+\lambda _{3}^{3}=-\alpha \delta -\delta ^{2},%
\end{array}
\\
&&%
\begin{array}{c}
c=-\alpha ^{2}-\delta ^{2},%
\end{array}
\\
&&%
\begin{array}{c}
\lambda _{1}^{2}=\lambda _{1}^{3}=0.%
\end{array}%
\end{eqnarray*}%
Hence, 
\begin{equation*}
\lambda _{2}^{2}=\delta \left( \delta -\alpha \right) \text{ and }\lambda
_{3}^{3}=\alpha \left( \alpha -\delta \right) .
\end{equation*}

In the remaining case, $\left( \beta ,\gamma \right) =\left( 0,0\right) $
and $\delta =\alpha $ and so, the algebraic Ricci Soliton condition 
\eqref{Soliton}
on $G_{5}$ is satisfied if and only if 
\begin{equation}
\lambda _{1}^{2}=\lambda _{1}^{3}=\lambda _{2}^{2}=\lambda _{2}^{3}=\lambda
_{3}^{2}=\lambda _{3}^{3}=0\text{ and }c=-2\alpha ^{2}. \square  \notag
\end{equation}

\begin{remark} 
Let $D$ denote the derivation of $\mathfrak{g}_{6}$ given in the above Theorem. 
Then it is easy to prove that $D$ is an outer derivation. 
\end{remark}

\subsection{Algebraic Ricci Solitons of $G_{7}$}
Finally consider the three-dimensional Lorentzian Lie group $G_{7}$. We
recall the following result of \cite{C07a}:

\begin{lemma}[\protect\cite{C07a}]
Let $\{e_{1},e_{2},e_{3}\}$ be the pseudo-orthonormal basis used in 
\eqref{g7}%
. Then%
\begin{eqnarray*}
\nabla _{e_{1}}e_{1} &=&\alpha e_{2}+\alpha e_{3},\text{ \ \ \ \ }\nabla
_{e_{2}}e_{1}=\beta e_{2}+\left( \beta +\frac{\gamma }{2}\right) e_{3},\text{
\ }\nabla _{e_{3}}e_{1}=-\left( \beta -\frac{\gamma }{2}\right) e_{2}-\beta
e_{3}, \\
\nabla _{e_{1}}e_{2} &=&-\alpha e_{1}+\frac{\gamma }{2}e_{3},\text{ \ }%
\nabla _{e_{2}}e_{2}=-\beta e_{1}+\delta e_{3},\text{ \ \ \ \ \ \ \ \ \ }%
\nabla _{e_{3}}e_{2}=\left( \beta -\frac{\gamma }{2}\right) e_{1}-\delta
e_{3}, \\
\nabla _{e_{1}}e_{3} &=&\alpha e_{1}+\frac{\gamma }{2}e_{2},\text{ \ \ \ }%
\nabla _{e_{2}}e_{3}=\left( \beta +\frac{\gamma }{2}\right) e_{1}+\delta
e_{2},\text{ \ \ }\nabla _{e_{3}}e_{3}=-\beta e_{1}-\delta e_{2},
\end{eqnarray*}%
and the only possibly non-zero components of the curvature tensor are given
by%
\begin{eqnarray*}
R_{1212} &=&\alpha \delta -\alpha ^{2}-\beta \gamma -\frac{\gamma ^{2}}{4},%
\text{ \ \ }R_{1313}=\alpha \delta -\alpha ^{2}-\beta \gamma +\frac{\gamma
^{2}}{4},\text{ \ } \\
\text{\ \ \ }R_{1213} &=&\alpha ^{2}-\alpha \delta +\beta \gamma ,\text{ \ \
\ \ \ \ \ \ \ }R_{2323}=-\frac{3}{4}\gamma ^{2}.
\end{eqnarray*}%
The Ricci operator is given by%
\begin{equation*}
\mathrm{Ric}=\left( 
\begin{array}{ccc}
-\frac{\gamma ^{2}}{2} & 0 & 0 \\ 
0 & \alpha \delta -\alpha ^{2}-\beta \gamma +\frac{\gamma ^{2}}{2} & \alpha
^{2}-\alpha \delta +\beta \gamma \\ 
0 & -\alpha ^{2}+\alpha \delta -\beta \gamma & -\alpha \delta +\alpha
^{2}+\beta \gamma +\frac{\gamma ^{2}}{2}%
\end{array}%
\right) .
\end{equation*}
\end{lemma}

Now, let $\mathrm{D}e_{l}=\lambda _{l}^{1}e_{1}+\lambda
_{l}^{2}e_{2}+\lambda _{l}^{3}e_{3}$ be an endomorphism of $\mathfrak{g}_{7}$
where $\{e_{1},e_{2},e_{3}\}$ is the pseudo-orthonormal basis used in 
\eqref{g7}%
. Then, $\mathrm{D}\in \mathrm{Der}\left( \mathfrak{g}_{7}\right) $ if and
only if%
\begin{eqnarray}
&&%
\begin{array}{c}
-\alpha \lambda _{2}^{2}+\alpha \lambda _{2}^{3}-\gamma \lambda
_{1}^{3}+\beta \lambda _{2}^{1}+\beta \lambda _{3}^{1}=0,%
\end{array}
\label{sysD7} \\
&&%
\begin{array}{c}
-\beta \lambda _{1}^{1}+\beta \lambda _{2}^{3}-\delta \lambda
_{1}^{3}+\alpha \lambda _{1}^{2}+\beta \lambda _{3}^{2}=0,%
\end{array}
\notag \\
&&%
\begin{array}{c}
\beta \left( -\lambda _{1}^{1}-\lambda _{2}^{2}+\lambda _{3}^{3}\right)
+2\beta \lambda _{2}^{3}+\left( \alpha -\delta \right) \lambda _{1}^{3}=0,%
\end{array}
\notag \\
&&%
\begin{array}{c}
\alpha \lambda _{3}^{3}-\alpha \lambda _{3}^{2}+\gamma \lambda
_{1}^{2}-\beta \lambda _{2}^{1}-\beta \lambda _{3}^{1}=0,%
\end{array}
\notag \\
&&%
\begin{array}{c}
\beta \left( \lambda _{1}^{1}-\lambda _{2}^{2}+\lambda _{3}^{3}\right)
-2\beta \lambda _{3}^{2}+\left( \delta -\alpha \right) \lambda _{1}^{2}=0,%
\end{array}
\notag \\
&&%
\begin{array}{c}
\beta \lambda _{1}^{1}-\beta \lambda _{3}^{2}+\delta \lambda _{1}^{2}-\alpha
\lambda _{1}^{3}-\beta \lambda _{2}^{3}=0,%
\end{array}
\notag \\
&&%
\begin{array}{c}
\gamma \left( -\lambda _{1}^{1}+\lambda _{2}^{2}+\lambda _{3}^{3}\right)
+\left( \alpha -\delta \right) \left( \lambda _{2}^{1}+\lambda
_{3}^{1}\right) =0,%
\end{array}
\notag \\
&&%
\begin{array}{c}
\delta \lambda _{3}^{3}+\beta \lambda _{2}^{1}+\beta \lambda _{3}^{1}-\gamma
\lambda _{1}^{2}-\delta \lambda _{3}^{2}=0,%
\end{array}
\notag \\
&&%
\begin{array}{c}
\delta \lambda _{2}^{2}+\beta \lambda _{2}^{1}+\beta \lambda _{3}^{1}-\gamma
\lambda _{1}^{3}-\delta \lambda _{2}^{3}=0.%
\end{array}
\notag
\end{eqnarray}%
In determining the solutions of 
\eqref{sysD7}
we must also take into account that, by 
\eqref{g7}%
, $\alpha +\delta \neq 0$ and $\alpha \gamma =0.$ Routine but long
calculation lead to the following:

\begin{lemma}
The solutions of 
\eqref{sysD7}
are the following:

\begin{itemize}
\item If $\delta \neq \alpha :$%
\begin{equation*}
\mathrm{D=}\left( 
\begin{array}{ccc}
\lambda _{1}^{1} & \lambda _{2}^{1} & \frac{\gamma }{\left( \alpha -\delta
\right) }\left( \lambda _{1}^{1}-\lambda _{2}^{2}-\lambda _{3}^{3}\right)
-\lambda _{2}^{1} \\ 
\frac{\beta }{\left( \alpha -\delta \right) }\left( \lambda _{1}^{1}-\lambda
_{2}^{2}-\lambda _{3}^{3}\right) & \lambda _{2}^{2} & \lambda _{3}^{3} \\ 
\frac{\beta }{\left( \alpha -\delta \right) }\left( \lambda _{1}^{1}-\lambda
_{2}^{2}-\lambda _{3}^{3}\right) & \lambda _{2}^{2} & \lambda _{3}^{3}%
\end{array}%
\right) .
\end{equation*}

\item If $\delta =\alpha \ $and $\beta \neq 0:$%
\begin{equation*}
\mathrm{D=}\left( 
\begin{array}{ccc}
\lambda _{1}^{1} & \lambda _{2}^{1} & -\lambda _{2}^{1} \\ 
\lambda _{1}^{2} & \lambda _{2}^{2} & \lambda _{1}^{1}-\lambda _{2}^{2} \\ 
\lambda _{1}^{2} & \lambda _{2}^{2} & \lambda _{1}^{1}-\lambda _{2}^{2}%
\end{array}%
\right) .
\end{equation*}

\item If $\delta =\alpha $ and $\beta =0:$%
\begin{equation*}
\mathrm{D=}\left( 
\begin{array}{ccc}
\lambda _{1}^{1} & \lambda _{2}^{1} & -\lambda _{2}^{1} \\ 
\lambda _{1}^{2} & \lambda _{2}^{2} & \lambda _{3}^{3} \\ 
\lambda _{1}^{2} & \lambda _{2}^{2} & \lambda _{3}^{3}%
\end{array}%
\right) .
\end{equation*}
\end{itemize}
\end{lemma}

Using this result, we now prove the following.

\begin{theorem}
Consider the three-dimensional Lorentzian Lie group $G_{7}$. Then, $G_{7}$
is an algebraic Ricci Soliton Lie group if and only if 
\begin{equation*}
\mathrm{D=}\alpha \left( \delta -\alpha \right) \left( 
\begin{array}{ccc}
0 & 0 & 0 \\ 
0 & 1 & -1 \\ 
0 & 1 & -1%
\end{array}%
\right), c=0 \text{ and } \gamma =0.
\end{equation*}
\end{theorem}

\textbf{Proof.} Assume that $\alpha \neq \delta $. Then the algebraic Ricci
Soliton condition 
\eqref{Soliton}
on $G_{7}$ is satisfied if and only if 
\begin{eqnarray}
&&%
\begin{array}{c}
c+\lambda _{2}^{2}=-\left( c+\lambda _{3}^{3}\right) =\alpha \left( \delta
-\alpha \right) -\beta \gamma +\frac{\gamma ^{2}}{2},%
\end{array}
\label{Sol7a} \\
&&%
\begin{array}{c}
\beta \left( \lambda _{1}^{1}-\lambda _{2}^{2}-\lambda _{3}^{3}\right) =0,%
\end{array}
\notag \\
&&%
\begin{array}{c}
\lambda _{2}^{2}=-\lambda _{3}^{3}=\alpha \left( \delta -\alpha \right)
-\beta \gamma +\frac{\gamma ^{2}}{2},%
\end{array}
\notag \\
&&%
\begin{array}{c}
\gamma \left( \lambda _{1}^{1}-\lambda _{2}^{2}-\lambda _{3}^{3}\right)
=\lambda _{2}^{1}=0,%
\end{array}
\notag \\
&&%
\begin{array}{c}
c+\lambda _{1}^{1}=-\frac{\gamma ^{2}}{2},%
\end{array}
\notag
\end{eqnarray}

Hence, it follows from the first and the third equation in 
\eqref{Sol7a}
that $c=0$ and so, $\lambda _{1}^{1}=-\frac{\gamma ^{2}}{2}.$ This means
that, where we took into account the fourth equation of 
\eqref{Sol7a}%
, 
\eqref{Sol7a}
holds if and only if $\gamma =0.$ Thus, we get 
\begin{equation*}
\lambda _{2}^{2}=-\lambda _{3}^{3}=\alpha \left( \delta -\alpha \right) 
\text{ and }\lambda _{2}^{1}=\lambda _{1}^{1}=c=0.
\end{equation*}

Now, suppose that $\delta =\alpha $ and $\beta $ $\neq 0.$ Note that in this
case $\gamma =0$ (since $\alpha +\delta \neq 0$ and $\alpha \gamma =0$) and
so \textrm{Ric}$=0$. Therefore, the algebraic Ricci Soliton condition 
\eqref{Soliton}
on $G_{7}$ is satisfied if and only if%
\begin{equation*}
\lambda _{1}^{1}=\lambda _{1}^{2}=\lambda _{2}^{1}=\lambda _{2}^{2}=c=0.
\end{equation*}

Next, assume that $\delta =\alpha $ and $\beta $ $=0.$ Thus, the algebraic
Ricci Soliton condition 
\eqref{Soliton}
on $G_{7}$ reduces to%
\begin{eqnarray}
&&%
\begin{array}{c}
c+\lambda _{1}^{1}=c+\lambda _{3}^{3}=-\left( c+\lambda _{2}^{2}\right) =-%
\frac{\gamma ^{2}}{2},%
\end{array}
\label{Sol7b} \\
&&%
\begin{array}{c}
\lambda _{1}^{2}=\lambda _{2}^{1}=\lambda _{2}^{2}=\lambda _{3}^{3}=0.%
\end{array}
\notag
\end{eqnarray}

Hence, 
\eqref{Sol7b}
holds if and only if $\gamma =0.$ Note that, the above two cases can be
deduced from the case when $\alpha \neq \delta $. $\square $


\begin{remark}
Let $D$ denote the derivation of $\mathfrak{g}_{7}$ given in the above
Theorem. Thus,


\begin{itemize}
\item If $\delta \neq 0$, $D$ satisfies 
\begin{equation*}
D=\dfrac{\alpha (\alpha -\delta )}{\delta }(\mathrm{ad}e_{2}+\mathrm{ad}%
e_{3}).
\end{equation*}

\item If $\delta =0$, in this case $D\neq \mathrm{ad}X$, for any $X\in 
\mathfrak{g}_{7}$.
\end{itemize}

\end{remark}











\end{document}